\documentclass[preprint,11pt]{elsarticle}

\usepackage[english]{babel}
\usepackage[utf8]{inputenc}
\usepackage[svgnames]{xcolor}
\usepackage[T1]{fontenc}

\usepackage{here}
\usepackage{color}
\usepackage{ragged2e}
\usepackage{array}
\usepackage{indentfirst}
\usepackage[hidelinks]{hyperref}
\usepackage{setspace}

\usepackage{amssymb}
\usepackage{amsmath}
\usepackage{amsthm}
\usepackage{amsfonts}

\usepackage{bm} 
\usepackage{stmaryrd} 
\usepackage{dsfont} 
\usepackage{mathrsfs}
\usepackage{enumitem} 
\usepackage{subcaption} 
\usepackage{multicol} 
\usepackage{multirow} 

\biboptions{sort}
\allowdisplaybreaks 

\newtheoremstyle{notation}
  {}
  {}
  {}
  {}
  {\bfseries}
  {}
  {7pt}
  {\thmname{#1}\thmnumber{ #2}.\textnormal{\thmnote{ (#3)}}}

\newtheoremstyle{exem}
  {}
  {}
  {}
  {}
  {\bfseries}
  {}
  { }
  {\textit{\thmname{#1}}\thmnumber{ #2}.\textnormal{\thmnote{ (#3)}}}

\theoremstyle{plain}
\newtheorem{theo}{Theorem}[section]
\newtheorem*{theo*}{Theorem}
\newtheorem{conj}[theo]{Conjecture}
\newtheorem{lemme}[theo]{Lemma}
\newtheorem*{lemme*}{Lemma}

\newtheorem*{sublemme*}{Sublemma}
\newtheorem{prop}[theo]{Proposition}
\newtheorem*{prop*}{Proposition}

\newtheorem{defi}[theo]{Definition}
\newtheorem*{defi*}{Definition}

\theoremstyle{definition}
\newtheorem{rem}[theo]{Remark}
\theoremstyle{notation}
\newtheorem*{nota}{Notation}
\theoremstyle{remark}

\numberwithin{equation}{section}
\numberwithin{figure}{section} 
\counterwithin{table}{section}

\newcommand{\pscal}[3][]{\left\langle #2,#3 \right\rangle_{#1}} 
\newcommand{\partent}[1]{\left\lfloor #1 \right\rfloor} 

\newcommand{\N}{\mathbb{N}} 
\newcommand{\Z}{\mathbb{Z}} 
\newcommand{\Q}{\mathbb{Q}} 
\newcommand{\R}{\mathbb{R}} 

\newcommand{\dimH}{\operatorname{dim}_H} 

\newcommand\numberthis{\addtocounter{equation}{1}\tag{\theequation}} 

\journal{Journal of Mathematical Analysis and Applications}

\begin{document}

\begin{frontmatter}



\title{Prevalent smoothness in inhomogeneous Besov spaces}


\author[labelLama]{Quentin Rible} 

\affiliation[labelLama]{organization={Univ Paris Est Creteil, Univ Gustave Eiffel, CNRS, LAMA UMR8050},
            postcode={F-94010 Creteil},
            country={France}}

\begin{abstract}
    In this article, we prove that, under some assumptions on the so-called environment, prevalent functions in inhomogeneous Besov spaces recently introduced by Barral-Seuret in 2023 are multifractal, with a singularity spectrum that we determine. This completes the previous Baire generic results already obtained.
\end{abstract}



\end{frontmatter}


\section{Introduction}\label{sec:intro}

The interest in multifractal properties dates back to the 1980s, when physicists sought to quantify the local variations in the velocity of a turbulent fluid. Indeed, turbulent flows are not spatially homogeneous, and the pointwise regularity of the velocity seems to differ widely from point to point. The relationship between multifractality and the global scaling properties of the velocity was first discovered by U. Frisch and G. Parisi \cite{Frisch-Parisi:1985:Multifractal-turbulence}. Since then, multifractal notions have proved to be useful and relevant in various mathematical domains: geometric measure theory, real analysis \cite{Jaffard:1996:Riemann_function_singularities}, dynamical systems, ergodic theory, diophantine approximation \cite{Feng:2007:Gibbs_properties_conformal_measures, Shmerkin:2005:multifractal_formalism_self-similar_measures}, to name a few. 

This article focuses on the multifractal analysis of functions, which consists of describing the pointwise behaviours of functions $f:\R^d\to\R$, computing their singularity spectrum $\sigma_f$, and, when possible, relating $\sigma_f$ with the scaling function associated with $f$. Given $x_0\in\R^d$, $H\in\R_+$, a function $f:\R^d\to\R$ belongs to $\mathcal{C}^{H}(x_0)$ if there exist a polynomial $P_{x_0}$ of degree less than $\partent{H}$, a constant $C>0$, and a neighbourhood $V$ of $x_0$ such that
\begin{equation}\label{eq:hold_spc}
    \forall x\in V, \quad |f(x)-P_{x_0}(x-x_0)| \leq C |x-x_0|^H.
\end{equation}
The \textit{pointwise Hölder exponent} of $f\in L^{\infty}_{\operatorname{loc}}(\R^d)$ at $x_0$ is then
\begin{equation}\label{eq:hold_exp}
    h_f(x_0)=\sup\{H\in\R_+ : f\in\mathcal{C}^H(x_0)\}.
\end{equation}
The level sets of $x\mapsto h_f(x)$ are $E_f(h) := \{x\in \R^d:h_f(x)=h\}$, and the associated \textit{singularity} (or \textit{multifractal}) \textit{spectrum} $\sigma_f$ of $f$ is the mapping
\begin{equation}\label{eq:defi:mf_spect}
    \sigma_f: h\in\R\cup\{\infty\} \mapsto \dimH(E_f(h)) \in [0,d] \cup \{-\infty\}.
\end{equation}
Here, $\dim_{H}$ stands for the Hausdorff dimension with $\dimH(\emptyset) =-\infty$.
This singularity spectrum describes the size of the level sets $E_f(h)$ of the singularities of $f$ in terms of Hausdorff dimension and provides a geometrical description of their distribution.

Results on the pointwise regularity of functions in classical spaces such as Bessel potential spaces $H^{s,p}(\R^d)$, Sobolev-Slobodeckij $W^{s,p}(\R^d)$ or Besov $B^{s}_{p,q}(\R^d)$ spaces have been obtained in \cite{Aubry-Bastin-Dispa:2007:prevalence_sobolev, Fraysse-Jaffard:2006:Almost_all_function_sobolev, Jaffard:2000:Frisch-Parisi-conjecture, Jaffard:1997:multifractal_formalism_all_functions, Jaffard-Meyer:2000:regularity_critical_besov_spaces}. They are of two types: first, the identification of a universal upper bound for the singularity spectrum of all functions in a given space, and second, when possible, the computation of the explicit value for the spectrum of a generic set of functions in the sense of prevalence or Baire categories. 

Let us emphasize that for large classes of functions $f$ (and measures), $\sigma_f$ satisfies a so-called \textit{multifractal formalism}, i.e. $\sigma_f$ coincides with the Legendre transform $\tau_{f}^*$ of a scaling function $\tau_{f}: \R \to \R$ associated with $f$,
\begin{equation}\label{eq:mf_form}
    \sigma_{f}(h) = \tau_{f}^*(h) := \inf_{q\in\R} hq - \tau_{f}(q).
\end{equation}
The scaling function $\tau_{f}$, defined using wavelets (see \cite{Jaffard:2004:Wavelet_techniques, Jaffard:2007:Wavelet_leaders}), describes the statistical scaling properties of $f$. This relationship is key, since when properly defined, $\tau_f$ is numerically estimable on real data, see \cite{Jaffard:2007:Wavelet_leaders}. Such estimation algorithms have generated breakthrough methods in signal and image processing \cite{Abry-Jaffard-Wendt:2009:wavelet_leader_texture_classification, Wendt-Abry-Jaffard-Roux:2009:Wavelet_Bootstrap} and other scientific fields \cite{Wendt-Abry-Ciuciu-Dumeur-Jaffard-Saes:2024:weak_scaling_MEG}. When \eqref{eq:mf_form} holds, one says that the \textit{multifractal formalism} holds for $f$.

A key feature of typical functions in the classical spaces $W^{s,p}(\R^d)$ and $B^{s}_{p,q}(\R^d)$ is that their singularity spectrum is always affine increasing. This is an apparent contradistinction with the fact that the spectra estimated on real data through the multifractal formalism always have a strictly concave shape, with an increasing {\em and} a decreasing part. 

Jaffard \cite{Jaffard:2000:Frisch-Parisi-conjecture} showed that given an {\em admissible} concave function $\tau:\R^+\to \R$, Baire typical functions in $\bigcap_{\varepsilon>0, p>0} B^{(\tau(p)-\varepsilon)/p}_{p,q,loc}(\R^d)$ have concave increasing spectra.
Fraysse \cite{Fraysse:2007:prevalent_validity_multifractal_formalism} obtained that prevalent functions in such spaces possess the same spectrum.

In order to consider non-affine singularity spectra, one could look at the existing natural generalisations of the classical Sobolev and Besov spaces, the weighted Sobolev and Besov spaces, which have various possible definitions.

For $h\in\R^d$ and $f:\R^d \to \R$, consider the finite difference operator $\Delta_{h}f:x\in\R^d\mapsto f(x+h)-f(x)$. Then, for $n\geq 2$, set $\Delta^{n}_{h}f=\Delta_{h}(\Delta^{n-1}_{h}f)$.
    
Let $\rho:\R\to\R$ be a measurable function, positive almost everywhere, called a weight. Under some hypotheses on the behaviour of $\rho$ given in \cite{Ansorena-Blasco:1995:weighted_besov_spaces}, weighted Besov spaces $B^{\rho}_{p,q}(\R^d)$ consist of functions $f:\R^d \to \mathbb{C}$ such that
\begin{equation*}
    \int_{|h| \leq 1} \rho(|h|)^{-q} \|\Delta_{h}^{n}f\|_{L^p(\R^d)}^q \frac{dh}{|h|^d} < \infty.
\end{equation*}
A function $f$ in such spaces presents a homogeneous pointwise behaviour as the weight in the integral depends only on $|h|$. 
These spaces have a generalised version called \textit{Besov space of generalised smoothness} $B^{\gamma,N}_{p,q}$, with $\gamma, N$ some sequences, studied by Farkas-Leopold \cite{Farkas-Leopold:2006:Characterisations_generalised_smoothness_Besov}. Moura \cite{Moura:2007:Characterisation_generalised_smoothness_besov} and Loosveldt-Nicolay \cite{Loosveldt-Nicolay:2019:equivalent_definition_generalised_smoothness_besov} gave alternative definitions. Kreit-Nicolay linked them to generalized Hölder-Zygmund $\Lambda^{\gamma}_{\sigma,N}$ spaces \cite{Kreit-Nicolay:2012:Characterisation_generalised_Holder_spaces}. 
These last spaces have been linked to traces of functions in Besov spaces on fractal sets by Caetano-Haroske \cite{Caetano-Haroske:2016:traces_besov_fractal} and their multifractal properties studied by Kleyntssens-Nicolay \cite{Kleyntssens-Nicolay:2024:Generalised_S_nu_spaces}.

A second version of weighted Besov spaces exists where the weight does not dictate the scale behaviour of the oscillations of the function $f$ but the local behaviour of $f$ itself and can be linked to standard Besov spaces, see \cite{Edmunds-Triebel:1996:function_spaces, Haroske-Triebel:1994:entropy_number_weighted_Besov_Triebel_1, Triebel:1983:theory_function_spaces_1}. It is defined as follows: for a weight $\rho\in L^{1}_{\operatorname{loc}}(\R)$, a function $f \in L^{1}_{\operatorname{loc}}(\R^d)$ belongs to $B^{s}_{p,q}(\R^d,\rho)$ if 
\begin{equation*}
    \int_{|h| \leq 1} |h|^{-sq} \|\rho~\Delta_{h}^{n}f\|_{L^p(\R^d)}^q \frac{dh}{|h|^d} < \infty.
\end{equation*}
Even if the multifractal behaviour of functions in such spaces has not been studied yet, they have an explicit characterisation through wavelets given by Triebel \cite[Chapter 6]{Triebel:2006:theory_function_spaces_3} that can be linked to the spaces introduced hereafter.

Another consideration is the variable exponent Lebesgue $L^{p(\cdot)}(\R^d)$ and Besov spaces $B^{s}_{p(\cdot),q(\cdot)}(\R^d)$, which are defined through the weighted Orlicz spaces \cite{Diening-Harjulehto-Hasto-Ruzicka:2011:variable_exponent_lebesgue_sobolev}, where the usual real parameters $p,q$ are replaced by continuous functions $p,q:\R^d \to [1,\infty)$. 
While typical functions in $B^{s}_{p,q}(\R^d)$ exhibit affine spectra with no decreasing parts \cite{Jaffard:2000:Frisch-Parisi-conjecture}, the spaces $B^{s}_{p(\cdot),q(\cdot)}(\R^d)$, have not been studied yet.

The functions in standard Besov and Sobolev spaces and the previous function spaces typically exhibit a concave singularity spectrum \cite{Jaffard:2000:Frisch-Parisi-conjecture, Fraysse-Jaffard:2006:Almost_all_function_sobolev} with \textit{only} an increasing part - even when it is not yet proved, there is not enough inhomogeneity in their definition to change this fact. To provide a natural framework in which typical functions have multifractal behaviours similar to those observed in real data, Barral and Seuret \cite{Barral-Seuret:2023:Besov_Space_Part_2} introduced \textit{inhomogeneous} Besov spaces, in which Baire typical functions may have as singularity spectrum any admissible continuous concave curve with an increasing and a decreasing part, prescribed in advance. 
For a H\"older set function $\mu$ (which can be thought as a measure for the moment, see Definition \ref{defi:hold_set_func} below), the \textit{inhomogeneous} Besov space $B^{\mu}_{p,q}(\R^d)$ consists of all functions $f$ in $L^p(\R^d)$ such that
\begin{equation}\label{eq:inh_besov_int}
    \int_{|h| \leq 1} |h|^{-dq/p} \bigg\|\frac{\Delta_{h}^{n}f}{\mu(B(\cdot+nh/2,nh/2))}\bigg\|_{L^p(\R^d)}^q \frac{dh}{|h|^d} < +\infty
\end{equation}

\begin{figure}[ht]
    \centering                                                                                                                   
    \includegraphics[page=5,width=.4\linewidth]{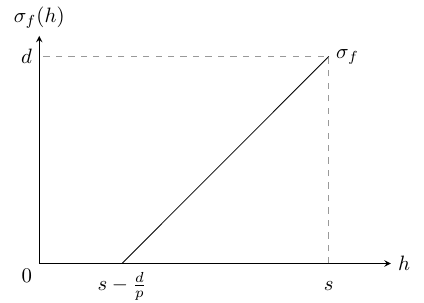}
    \includegraphics[page=2,width=.54\linewidth]{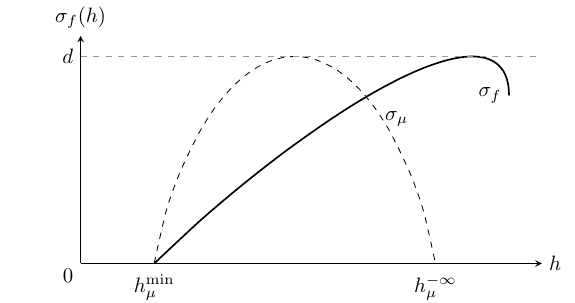}
    \caption{This figure illustrates the typical singularity spectrum in standard Besov spaces, as opposed to the typical one in inhomogeneous Besov spaces, to motivate their consideration. \textbf{Left:} Singularity spectrum of a typical function $f\in B^{s}_{p,q}(\R^d)$. \textbf{Right:} Singularity spectrum of a typical function $f\in B^{\mu}_{1,q}(\R^d)$. The dashed graph represents the (initial) singularity spectrum of $\mu$. When $p = \infty$ and $f$ is typical in $B^{\mu}_{\infty,q}(\R^d)$, $\sigma_f=\sigma_{\mu}$.}
    \label{fig:typic_spect_inh_besov}
\end{figure}
These spaces are not spatially homogeneous in the sense that, in \eqref{eq:inh_besov_int}, the oscillations $\Delta_{h}^{n}f (x)$ of $f$ are locally compared to the $\mu$-mass of a small ball of radius $h$ close to $x$, and this mass may strongly depend on $x$. 

In the present article, we obtain prevalent generic results for the multifractal properties of functions in inhomogeneous Besov spaces.

\medskip

We start with some definitions. The diameter of a set $E$ is written $|E|$.
\begin{defi}\label{defi:hold_set_func}
    The set of Hölder set functions is 
    \begin{equation*}
        \mathcal{H}(\R^d) = \{\mu:\mathcal{B}(\R^d)\to \R_+ \cup \{+\infty\} : \exists C, s > 0,\forall E \in \mathcal{B}(\R^d), \mu(E) \leq C |E|^s\}.
    \end{equation*}
\end{defi}

When $\mu(E) \leq C |E|^s $ for all $E$ with $0\leq |E| \leq 1$, $\mu$ is said to be $s$-H\"older.

The topological support of $\mu \in \mathcal{H}(\R^d)$ is $\operatorname{supp}(\mu)=\{x\in \R^d : \forall {r > 0},\break \mu(B(x, r)) > 0 \}$. 
A set function $\mu$ is fully supported on a Borel set $A$ when $\operatorname{supp}(\mu)=A$.
We will denote $\mathcal{H}([0,1]^d) = \{\mu \in \mathcal{H}(\R^d) : \operatorname{supp}(\mu) \subset [0,1]^d\}$.

For $\mu\in\mathcal{H}(\R^d)$, $s>0$, and $E\subset \R^d$, define $\mu^{s}(E):=\mu(E)^s$ and $\mu^{(+s)}(E):=\mu(E)|E|^s$,
and if $\mu$ is $s_0$-H\"older, define for all $s\in(0,s_0)$, 
\begin{equation*}
    \mu^{(-s)}(E):=
    \begin{cases}
        0               &\text{if}~|E|=0, \\
        \mu(E)|E|^{-s} &\text{if}~0<|E|<\infty.
    \end{cases}
\end{equation*}
Observe that $\mu^{s}$, $\mu^{(+s)}$ and $\mu^{(-s)}$, as defined above, still belongs to $\mathcal{H}(\R^d)$.

In this paper, we will need the additional assumption that the set functions $\mu\in\mathcal{H}(\R^d)$ are decreasing with respect to the embedding relationship, i.e. they are {\em capacities}. This additional property is natural in the functional setting. Indeed, looking back to \eqref{eq:inh_besov_int}, it is intuitive that if the oscillation of $f$ on an interval $I$ is controlled by $\mu(I)$, then the oscillations on a subinterval $J\subset I$ shall be controlled by $\mu(J) \leq \mu(I)$, since the oscillation on $J$ is necessarily smaller than that on $I$.
Hence, we define the sets of capacities
\begin{align*}
    \mathcal{C}(\R^d) &= \{\mu \in \mathcal{H}(\R^d) : \forall E,F \in \mathcal{B}(\R^d), E\subset F \Rightarrow \mu(E)\leq \mu(F)\}\\
    \mathcal{C}([0,1]^d) &= \mathcal{H}([0,1]^d) \cap \mathcal{C}(\R^d).
\end{align*}

\begin{defi}\label{defi:P1_ball}
	A set function $\mu\in\mathcal{H}(\R^d)$ satisfies property $(P)$ if there exist $C, s_1, s_2 >0$ such that
	\begin{equation}\label{eq:P1_ball}
		\text{for all $r>0$, for all $x\in\operatorname{supp}(\mu)$,}~~C^{-1} r^{s_2} \leq \mu(B(x,r)) \leq C r^{ s_1}.
	\end{equation}
\end{defi}

In particular, such a $\mu$ is $s_1$-H\"older.
We will only consider set functions satisfying property $(P)$ on $\R^d$ (or $[0,1]^d$), implying that these functions will be fully supported on $\R^d$ (or $[0,1]^d$). We will consider a capacity $\mu$ said to be \textit{almost doubling}.
Recall that a finite Borel measure $\mu$ is said to be doubling if for some constant $C_{\mu}\geq 1$, every $x\in \R^d$, $r>0$, 
\begin{equation*}
    \mu(B(x,r)) \leq \mu(B(x,2r))\leq C_\mu\,\mu(B(x,r)).
\end{equation*}

\begin{defi}\label{defi:alm_doub_rad}
    Let $\Phi$ be the set of non increasing functions $\phi:(0,1]\to \R_+$ such that $\lim_{r\to 0^+} \frac{\phi(r)}{\log(r)}=0$.
    A capacity $\mu\in\mathcal{C}(\R^d)$ is \textnormal{almost doubling} on a set $A$ when there exists $\phi\in\Phi$ such that for every $x\in A$, $r\in(0,1]$,
    \begin{equation}\label{eq:alm_doub_rad}
        \mu(B(x,2r)) \leq e^{\phi(r)} \mu(B(x,r)).
    \end{equation}
    The set function $\mu$ is doubling when \eqref{eq:alm_doub_rad} holds with a constant function $\phi$.
\end{defi}

Intuitively, \eqref{eq:alm_doub_rad} implies that moving the centre and radius of a ball slightly keeps its $\mu$-mass controlled. A doubling measure satisfies \eqref{eq:alm_doub_rad} with a constant function $\phi$. This almost doubling property generalises the concept of doubling measure while not impacting the multifractal behaviour of such set functions.
Examples of doubling set functions include (1-periodic) Gibbs measures $\mu$ \cite{Bowen:2008:ergodic_theory, Pesin-Weiss:1997:multi-frac-anal-gibbs-measures, Collet-Lebowitz-Porzio:1987:dimension_spectrum_dynamical_system} and their powers $\mu^s$ for $s>0$.

In \cite{Barral-Seuret:2023:Besov_Space_Part_2}, Barral and Seuret identified the typical multifractal behaviour, in the Baire category sense, of functions in $B^{\mu}_{p,q}(\R^d)$, when $\mu$ belongs to a large class of set functions satisfying $(P)$ but also an extra condition, see \cite{Barral-Seuret:2023:Besov_Space_Part_2}. 
Barral and Seuret used this space to solve the so-called Frisch-Parisi conjecture, which consists of proving the existence of function spaces in which typical functions satisfy the multifractal formalism \cite{Barral-Seuret:2023:Besov_Space_Part_2, Jaffard:1997:multifractal_formalism_all_functions} for any possible admissible spectrum.
To achieve this, they used a wavelet characterisation of functions, as we will in this article.

Write $L^d=\{0,1\}^d\setminus\{0^d\}$, and consider a wavelet family $\{\Phi,\Psi_\lambda:=\Psi^l(2^j \cdot -k)\}_{ \lambda=(j,k,l) \in \Z\times \Z^d\times L^d }$ generated by a family of wavelets $(\Phi,\Psi^l)_{l\in L^d}$ (see Section \ref{sec:mr_wlt_ana} for details), such that every $f\in L^2(\R^d)$ can be written
\begin{equation}\label{eq:wlt_dec_in_L2}
    f \overset{L^2}{=} \sum_{\lambda \in \Z\times \Z^d\times L^d} c_{\lambda} \Psi_{\lambda} \quad \text{ with }\quad c_{\lambda}=2^{dj} \int_{\R^d} f(x)\Psi_{\lambda} (x)dx.
\end{equation}
For $k=(k_1,\ldots,k_d)\in\Z^d$, $j\in\Z$, we often identify the triplet $\lambda=(j,k,l)$ with the dyadic cube $\prod_{i=1}^{d} [k_i 2^{-j}, (k_i+1) 2^{-j})$, and in particular we write $\mu(\lambda)$ for $\mu(\prod_{i=1}^{d} [k_i 2^{-j}, (k_i+1) 2^{-j}))$.

Write $\Lambda_j= \{ \lambda=(j,k,l) : k\in \Z^d, l\in L^d \}$ for the set of the dyadic cubes of generation $j$.
If the function $f$ has a support in $[0,1]^d$, we will write
\begin{equation}\label{eq:wlt_dec_in_01}
    f = \pscal{f}{\Phi} \, \Phi + \sum_{j\geq 0} \sum _{\lambda\in \Lambda_j} c_{\lambda} \Psi_{\lambda} .
\end{equation}

\begin{defi}\label{defi:inh_besov_wlt}
    Let $\mu\in\mathcal{H}(\R^d)$ with $\operatorname{supp}(\mu)=\R^d$ satisfying $(P)$ with $0<s_1<s_2$ and $p\in[1,+\infty]$. Consider an integer $n \geq \partent{s_2+\frac{d}{p}}+1$.
    Consider wavelets $(\Phi,\Psi^l)_{ l\in L^d}$ that belong to $C^n(\R^d)$ and have at least $n$ vanishing moments (see Section \ref{sec:mr_wlt_ana} for details).

    A function $f$ belongs to $B^{\mu,\Psi}_{p,q}(\R^d)$ if $\|f\|_{L^p(\R^d)} + |f|_{\mu,p,q,\Psi}<+\infty$, where 
    \begin{equation}\label{eq:snorm_inh_besov}
        |f|_{\mu,p,q,\Psi}=\big\|(\varepsilon^{\mu,p}_{j})_{j\in\N}\big\|_{\ell^q(\N)} \text{with}~ \varepsilon^{\mu,p}_{j}=\bigg\|\Big( \frac{c_{\lambda}}{\mu(\lambda)} \Big)_{\lambda\in \Lambda_j \times L^d}\bigg\|_{\ell^p(\Lambda_j \times L^d)}.
    \end{equation}
    
    Then define 
    \begin{equation}\label{eq:besov_tilde_wlt}
        \widetilde{B}^{\mu,\Psi}_{p,q}(\R^d) = \bigcap_{0 < \varepsilon < \min(1,s_1)} B^{\mu^{(-\varepsilon)},\Psi}_{p,q}(\R^d).
    \end{equation}
    
    Similarly, when $\operatorname{supp}(\mu)=[0,1]^d$, a function $f$ belongs to $B^{\mu,\Psi}_{p,q}([0,1]^d)$ if $\|f\|_{L^p(\R^d)} + |f|_{\mu,p,q,\Psi}<+\infty$, with the convention that $\frac{c_{\lambda}}{\mu(\lambda)} =0$ as soon as $\mu(\lambda)=0$, and $\widetilde{B}^{\mu,\Psi}_{p,q}([0,1]^d)$ is defined by replacing $\R^d$ by $[0,1]^d$ in \eqref{eq:besov_tilde_wlt}.
\end{defi}
Remark that $B^{\mu,\Psi}_{p,q}(\R^d) \subset \widetilde{B}^{\mu,\Psi}_{p,q}(\R^d)$ and similarly on $[0,1]^d$. In particular, when $p=q=\infty$, $B^{\mu,\Psi}_{\infty,\infty}(\R^d) $ is the space of those functions $f$ for which 
$|c_\lambda | \leq C \,\mu(\lambda)$, so the space $B^{\mu,\Psi}_{\infty,\infty}(\R^d)$ can be understood as a inhomogeneous generalisation of H\"older $C^\alpha(\R^d)$ spaces, which are known to be characterised by the fact that $|c_\lambda | \leq C\,|\lambda|^\alpha$ for every $\lambda$. 
This justifies, once again, the name "environment" used for $\mu$, since the set function may have various local behaviours and the constraint $|c_{\lambda}| \leq C\,\mu(\lambda)$ is highly dependent on $\mu$ and the location of the cube $\lambda$.

When $p<\infty$, the space $B^{\mu,\Psi}_{p,\infty}(\R^d)$ can be seen as a inhomogeneous generalisation of Cald\'eron and Zygmund $T_{p}^\alpha(\R^d)$ spaces.
Let us also mention that the spaces $C^\alpha(\R^d)$ and $T_{p}^\alpha(\R^d)$ have been generalised by Kreit, Loosveldt and Nicolay \cite{Kreit-Nicolay:2012:Characterisation_generalised_Holder_spaces, Kreit-Nicolay:2018:Generalised_pointwise_Holder_sequences, Loosveldt-Nicolay:2022:Prevalent_pointwise_regularity_TPalpha} with the spaces $C^{\gamma}(\R^d)$ and $T^{\gamma}_{p,q}(\R^d)$ with $\gamma$ a sequence. The space $C^{\gamma}(\R^d)$ is characterised by $|c_{\lambda_j}|\leq C\gamma_j$ for every $\lambda_j\in\Lambda_j$ (for a sequence $\gamma=(\gamma_j)$ verifying, for all $j$, $C^{-1}\gamma_j\leq \gamma_{j+1} \leq C \gamma_j$.)

\medskip

Since multifractal properties are local, we will focus on behaviours of functions on the cube $[0,1]^d$, and will consider functions written under the form \eqref{eq:wlt_dec_in_01}. This is not a restriction, since any function $f$ written as in \eqref{eq:wlt_dec_in_L2} can be decomposed as a (countable) sum $f = \widetilde{f} + \sum_{\ell\in \Z^d} f_\ell(.+\ell)$,
where $\widetilde{f}$ is regular (corresponding to the low frequency part of $f$) and each $f_\ell$ has a decomposition of the form \eqref{eq:wlt_dec_in_01} on $[0,1]^d$. Then, each function $f_\ell(.+\ell)$ has an interesting behaviour only on the cube $\ell+[0,1]^d$, and the multifractal properties of $f$ will be deduced from the juxtaposition of those of the $f_\ell$ using that $\sigma_{f} = \sup_{\ell\in \Z^d} \sigma_{f_{\ell}}$.

\smallskip

{\bf From now on, we focus on set functions $\mu\in \mathcal{H}([0,1]^d)$ and on functions belonging to $B^{\mu,\Psi}_{p,q}([0,1]^d)$ and $\widetilde{B}^{\mu,\Psi}_{p,q}([0,1]^d)$ written as \eqref{eq:wlt_dec_in_01}.}

\smallskip

Our goal is to obtain the prevalent multifractal behaviour in $B^{\mu,\Psi}_{\infty,q}(\R^d)$ in larger classes of capacities than in \cite{Barral-Seuret:2023:Besov_Space_Part_2}. 
Recall that the prevalence theory, proposed by Christensen \cite{Christensen:1972:sets_Haar_measure_zero} and Hunt, Sauer and York \cite{Hunt-Sauer-Yorke:1992:prevalence} independently, supersedes the notion of Lebesgue for universally measurable sets in any real or complex topological vector space $E$ endowed with its Borel $\sigma$-algebra $\mathcal{B}(E)$. All Borel measures $\vartheta$ on $(E, \mathcal{B}(E))$ will be automatically completed, i.e. for every set $A\subset E$ and $B \in \mathcal{B}(E)$ such that the symmetric difference $A\Delta B$ is included in some $D\in\mathcal{B}(E)$ with $\vartheta(D) = 0$, one sets $\vartheta(A) := \vartheta(B)$.
\begin{defi}
    A set $A\subset E$ is universally measurable if it is measurable for any (completed) Borel measure on $E$.
    A universally measurable set $A\subset E$ is called Haar-null, or shy, if there exists a finite completed measure $\vartheta$ that is positive on some compact subset $K$ of $E$ and such that 
    \begin{equation}
        \text{for every } x\in E, \quad \vartheta(A+x)=0.
    \end{equation}
    The measure $\vartheta$ used is called \textit{transverse}. 
    A set that is included in a Haar-null universally measurable set is also called Haar-null, and the complement in $E$ of a Haar-null set $A$, denoted $A^{\complement}$, is called prevalent. 
\end{defi}

Hunt pioneered the use of prevalence in function spaces \cite{Hunt:1994:prevalence_continuous_functions}. 
Bayart-Heurteaux \cite{Bayart-Heurteaux:2013:Hausdorff_Dimension_Graphs_Prevalent_Continuous_Functions} and Balka-Darji-Elekes \cite{Balka-Darji-Elekes:2016:dimension_graph_continous_maps} studied the dimension of images of continuous functions.
Fraysse-Jaffard-Kahane \cite{Fraysse-Jaffard-Kahane:2005:propriete_analyse, Fraysse-Jaffard:2006:Almost_all_function_sobolev} and Aubry-Maman-Seuret \cite{Aubry-Maman-Seuret:2013:Traces_besov_results} studied, respectively, the prevalent regularity of functions and traces of functions in Besov spaces. 
Olsen studied the multifractal and $L^q$-dimension of prevalent measures on compact sets of $\R^d$ \cite{Olsen:2010:Lq_dimension_prevalent_measure, Olsen:2010:multifractal_dimension_prevalent_measure}.

Let us give the last definitions needed to state our second main theorem.

\begin{defi}\label{defi:mf_form}
    Let $\mu\in\mathcal{H}(\R^d)$. For $x\in \operatorname{supp}(\mu)$, the \textit{lower local dimension} of $\mu$ at $x$ is defined by
    \begin{equation}\label{eq:h_mu_ball}
        \underline{h}_{\mu}(x) = \liminf_{r \to 0^+} \frac{\log \mu(B(x,r))}{\log(r)}.
    \end{equation}
    The \textit{singularity spectrum} of $\mu$ is then the Hausdorff dimension of the level sets of the local dimension
    \begin{equation*}
        \sigma_{\mu}\ : \ h\in\R \mapsto \dimH(\underline{E}_{\mu}(h)) ~\text{where}~ \underline{E}_{\mu}(h)=\{x\in\operatorname{supp}(\mu) \ : \ \underline{h}_{\mu}(x)=h\}.
    \end{equation*}
    We also defined $E_{\mu}(h)$ as the set of points $x$ where the local dimension \eqref{eq:h_mu_ball} is a limit and equals $h$.
     
    Finally, the scaling function of $\mu\in\mathcal{H}([0,1]^d)$ is defined by
    \begin{equation}\label{eq:scaling_function_tau_mu}
        \tau_\mu(q)=\liminf_{j\to+\infty} \frac{1}{-j}\log_2 \sum_{\substack{\lambda\in\Lambda_j \cap [0,1]^d}} \mu(\lambda)^q.
    \end{equation}
\end{defi}

The convention $0^q=0$ is adopted. When $\operatorname{supp}(\mu)=[0,1]^d$, the singularity spectrum of $\mu$ is bounded above by the Legendre transform $\tau_\mu^*$ of $\tau_\mu$ (recall formula \eqref{eq:mf_form}) : for every $h\geq 0$, $\sigma_{\mu}(h) \leq \tau_{\mu}^*(h)$, see \cite{Brown_Michon_Peyriere:1992:Multifractal-analysis-measures,Jaffard:2000:Frisch-Parisi-conjecture,LevyVehel-Vojak:1998:Choquet-capacities}. 

\begin{defi}
    A set function $\mu\in\mathcal{H}([0,1]^d)$ satisfies the strong multifractal formalism (SMF) if for every $h\geq 0$, $\sigma_\mu(h)=\tau_\mu^*(h)=\dimH(E_{\mu}(h))$.
\end{defi}


\noindent \textbf{We consider the case $p=\infty$}, the case $p<\infty$ would require a finer study.

\begin{theo}\label{theo:prev_spect_inh_besov_inf_q}
    Let $\mu\in\mathcal{C}([0,1]^d)$ be an almost doubling capacity satisfying $(P)$ and verifying the SMF.
    For all $\displaystyle f\in \widetilde{B}^{\mu,\Psi}_{\infty,q}([0,1]^d)$,
    \begin{equation}
        \sigma_f(h) \leq 
        \begin{cases}
           \, \sigma_{\mu}(h) &\text{if}~h \leq \tau'_{\mu}(0) \\
         \ \ d &\text{if}~h > \tau'_{\mu}(0).
        \end{cases}
    \end{equation}
There is a prevalent set of functions $f\in \widetilde{B}^{\mu,\Psi}_{\infty,q}([0,1]^d)$ satisfying $\sigma_f=\sigma_{\mu}.$
\end{theo}

\vspace{-5mm}\begin{figure}[H]
    \centering
    \includegraphics[page=1,width=0.6\linewidth]{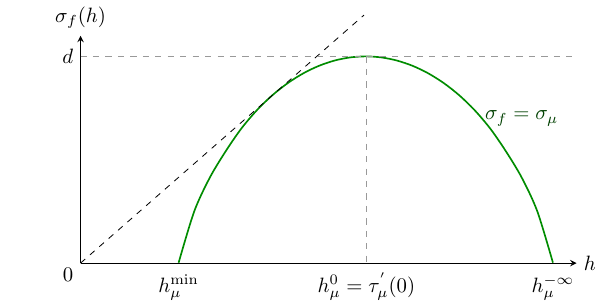}
    \caption{Example of singularity spectrum of a prevalent set of functions in $\widetilde{B}^{\mu,\Psi}_{\infty,q}([0,1]^d)$ for a measure $\mu$ with known singularity spectrum.}
    \label{fig:prevalent_spectrum}
\end{figure}

Remark that this result stays true for $f\in B^{\mu,\Psi}_{\infty,q}([0,1]^d)\subset \widetilde{B}^{\mu,\Psi}_{\infty,q}([0,1]^d)$. In other words, the prevalent singularity spectrum in $\widetilde{B}^{\mu,\Psi}_{\infty,q}$ or $B^{\mu,\Psi}_{\infty,q}$ coincides with that of $\mu$. 

Intuitively, almost all functions in $\widetilde{B}^{\mu,\Psi}_{\infty,q}$ or $B^{\mu,\Psi}_{\infty,q}$ have a concave increasing and decreasing singularity spectrum which is that of $\mu$. In \cite{Barral-Seuret:2023:Besov_Space_Part_2}, Barral and Seuret obtained the same singularity spectrum as in Theorem \ref{theo:prev_spect_inh_besov_inf_q} for Baire typical functions, but only for two large classes of capacities (a class of measures with a multifractal spectrum prescribed in advance and Gibbs capacities defined as powers of Gibbs measures). In this article, we work with a somewhat larger class of capacities (including the Gibbs capacities), which allows us to complete Barral-Seuret's result by stating that prevalent and Baire typical regularity coincide in this situation. As a consequence of our results, there exist wide classes of functions having, in the prevalent sense, a concave singularity spectrum having both an increasing and a decreasing part. Such behaviours fit those observed on real data, so these function spaces are natural frameworks for many physical phenomena exhibiting multifractal complex properties.

We shall prove Theorem \ref{theo:prev_spect_inh_besov_inf_q} for $\widetilde{B}^{\mu,\Psi}_{\infty,q}([0,1]^d)$, the proof for $B^{\mu,\Psi}_{\infty,q}([0,1]^d)$ is simpler, we will explain where to adapt our proof to get the simpler case.

It is certainly interesting to extend Theorem \ref{theo:prev_spect_inh_besov_inf_q} to the case $p<+\infty$.

\subsection*{Organisation of the paper}

Section \ref{sec:capa_prop_and_wlt_form} gives the necessary definitions and the properties of measures and capacities used to define a \textit{multifractal environment}, then presents the wavelet multifractal formalism and the wavelet leaders formalism developed by Jaffard in \cite{Jaffard:2004:Wavelet_techniques, Jaffard:2007:Wavelet_leaders}. 
In Section \ref{sec:Prevalent_spectrum_inf_q}, the proof of Theorem \ref{theo:prev_spect_inh_besov_inf_q} is given. 
For this, using a wavelet construction, we introduce a \textit{saturating function} $\mathscr{G}_q$ of $\widetilde{B}^{\mu,\Psi}_{\infty,q}([0,1]^{d})$ having a prescribed multifractal behaviour coinciding with that of $\mu$. This function $\mathscr{G}_q$ is a key object we build to show prevalent results on the spectrum of functions of $\widetilde{B}^{\mu,\Psi}_{\infty,q}([0,1]^{d})$.

\section{Additional definitions and properties on capacities, wavelets and multifractal formalism}\label{sec:capa_prop_and_wlt_form}

Let $\Lambda=\bigcup_{j\in\Z} \Lambda_j$ be the set of all dyadic cubes.

For $x\in\R^d$, $\lambda_j(x)$ is the dyadic cubes $\lambda$ of generation $j$ such that $x \in \lambda$, with the abuse of notation identifying $\lambda=(j,k=(k_1,...,k_d),l) \in \Z\times\Z^d\times L^d$ with the corresponding dyadic cube $\prod_{i=1}^{d} [k_i 2^{-j}, (k_i+1) 2^{-j})$.

For $j\in\Z$, $\lambda\in \Lambda_j$ and $N\in\N^*$, $N\lambda$ stands for the cube with same centre as $\lambda$ and radius equal to $N 2^{-j-1}$ in $(\R^d,\|\cdot\|_{\infty})$. For instance, $3\lambda$ is the union of those $\lambda'\in \Lambda_j$ such that $\partial\lambda \cap \partial \lambda' \neq \emptyset$ where $\partial \lambda$ is the boundary of $\lambda$.

\subsection{Capacity and formalism}\label{subsec:capa_and_form}

For any Borel set $A\subset\R^d$ and capacity $\mu$, we will denote $\mu\vert_{A} = \mu(.\cap A)$.
To describe the local behaviour of $\mu\in\mathcal{H}([0,1]^d)$, one needs, in addition to the lower local dimension \eqref{eq:h_mu_ball}, another local dimension.

\begin{defi}\label{defi:local_dim_and_set}
    Let $\mu\in\mathcal{H}([0,1]^d)$. For $x\in \operatorname{supp}(\mu)$, the upper local dimension of $\mu$ at $x$ is
    \begin{equation*}\label{eq:h_mu_sup_ball}
        \overline{h}_{\mu}(x)= \limsup_{r \to 0^+} \frac{\log \mu(B(x,r))}{\log(r)}.
    \end{equation*}
    Whenever $\underline{h}_{\mu}(x) = \overline{h}_{\mu}(x)$, the common limit is denoted by $h_{\mu}(x)$. For $\alpha\in\R$, 
    \begin{equation*}
    \begin{array}{c}
        \underline{E}_{\mu}(\alpha)=\{x\in\operatorname{supp}(\mu) : \underline{h}_{\mu}(x)=\alpha\},~ \overline{E}_{\mu}(\alpha)=\{x\in\operatorname{supp}(\mu) : \overline{h}_{\mu}(x)=\alpha\} \\ 
        \text{and} \quad E_{\mu}(\alpha) = \underline{E}_{\mu}(\alpha)\cap \overline{E}_{\mu}(\alpha).
    \end{array}
    \end{equation*}
\end{defi}

\begin{rem}\label{rem:local_dim_alt}
    Note that the local dimensions are sometimes defined as 
    \begin{align*}
        &\underline{h}_{\mu}(x)= \liminf_{j\to +\infty} \frac{\log_2 \mu(\lambda_j(x))}{-j} 
        ~\text{and}~ \overline{h}_{\mu}(x)= \limsup_{j\to +\infty} \frac{\log_2 \mu(\lambda_j(x))}{-j}, \numberthis \label{eq:h_mu_dyadic} \\
        \text{or}\quad 
        &\underline{h}_{\mu}(x)= \liminf_{j\to +\infty} \frac{\log_2 \mu(3\lambda_j(x))}{-j} 
        ~\text{and}~ \overline{h}_{\mu}(x)= \limsup_{j\to +\infty} \frac{\log_2 \mu(3\lambda_j(x))}{-j}.
    \end{align*}
    However, this paper mainly considers doubling or almost doubling set functions for which all previous notions of local dimensions \eqref{eq:h_mu_ball} and \eqref{eq:h_mu_sup_ball}, level sets, singularity spectrum and scaling function do not depend on whether dyadic cubes or centred balls are considered.
\end{rem}

Recall that the Legendre spectrum $h \mapsto \tau_{\mu}^{*}(h)$ is increasing on the interval $h\leq \tau_{\mu}'(0^{+})$ and is decreasing on $h\geq \tau_{\mu}'(0^{-})$. The following Proposition \ref{prop:prop_capa_SMF} can be deduced from \cite{Barral:2015:inverse-problems,Barral-Ben_Nasr-Peyriere:2003:multifract-formalism,Barral-Seuret:2023:Besov_Space_Part_2,LevyVehel-Vojak:1998:Choquet-capacities,Olsen:1995:Multifractal-formalism-measure}. 

\begin{prop}\label{prop:prop_capa_SMF}
    Let $\mu\in\mathcal{H}([0,1]^d)$ with $supp(\mu) = [0,1]^d$, such that $\mu$ obeys the SMF (Definition \ref{defi:mf_form}). Then for every $h\in\R^{+}$:
    \begin{enumerate}[label=$({\roman*})$]
        \item\label{item:prop_capa_SMF:sigma_is_tau_star} One has $\sigma_{\mu}(h) = \dimH(\underline{E}_{\mu}(h)) = \dimH(\overline{E}_{\mu}(h)) = \tau_{\mu}^{*}(h)$.
        \item\label{item:prop_capa_SMF:dim_E_leq} When $h\leq \tau'_{\mu}(0^{+})$, $\dimH(\underline{E}^{\leq}_{\mu}(h)) = \tau_{\mu}^{*}(h)$.
        \item\label{item:prop_capa_SMF:dim_E_geq} When $h\geq \tau'_{\mu}(0^{-})$, $\dimH(\overline{E}^{\geq}_{\mu}(h)) = \tau_{\mu}^{*}(h)$.
        \item\label{item:prop_capa_SMF:sigma_is_d} When $\tau'_{\mu}(0^{+}) \leq h \leq \tau'_{\mu}(0^{-})$, $\sigma_{\mu}(h) = d$.
    \end{enumerate}
\end{prop}

Set functions obeying such formalism are common (Gibbs measures \cite{Pesin-Weiss:1997:multi-frac-anal-gibbs-measures}, Mandelbrot cascades, self-similar measures associated with IFS satisfying open set conditions \cite{Brown_Michon_Peyriere:1992:Multifractal-analysis-measures, Patzschke:1997:self-conformal-measures}, ...), and constitute natural environments for inhomogeneous spaces.
Note that measures obeying the SMF are not almost doubling in general. Binomial measures obey the SMF \cite{Brown_Michon_Peyriere:1992:Multifractal-analysis-measures} but are not almost doubling, whereas Gibbs measures are almost doubling.

All the set functions we work with also satisfy the property $(P)$ given in Definition \ref{defi:P1_ball}. 
Observe that $(P)$ implies that for all $j\in\N$ and $\lambda\in\Lambda_j$,
\begin{equation}\label{eq:P1}
    C^{-1} 2^{-j s_2} \leq \mu(\lambda) \leq C 2^{-j s_1}.
\end{equation}
Property $(P)$ and \eqref{eq:P1} are not equivalent in general, but it is the case when $\mu$ is almost doubling (with slightly different constant $C$ and exponents $s_1$ and $s_2$), which will always be the case for us.
While Baire-typical measures on $[0,1]^d$ verify a multifractal formalism \cite{Buczolich-Seuret:2010:typical_measure_mf}, they do not satisfy $(P)$ in general. Some measures, such as sums of Dirac masses, do not verify either a multifractal formalism or $(P)$.

\subsection{Multi-resolution wavelet analysis}\label{sec:mr_wlt_ana}

Let $d\in\N^{*}$ and $1\leq p, q \leq +\infty$. A mapping $\psi:\R^d\to\R$ has $n+1$ vanishing moments when for every multi-index $\alpha\in\N^d$ such that $\sum_{i=1}^{d} \alpha_i \leq n$, $\int_{\R^d} x_1^{\alpha_1}\ldots x_d^{\alpha_d} \psi(x)dx=0$. 

For an arbitrary integer $N\geq 1$, there exist compactly supported functions $\Psi^{0} \in C^{N}(\R)$ (called the scaling function) and $\Psi^{1} \in C^{N}(\R)$ (called the mother wavelet), with $\Psi^{1}$ having at least $N+1$ vanishing moments and such that $(\Psi^{1}_{j,k}:x\mapsto \Psi^{1}(2^j x-k))_{j,k\in\Z}$ form an orthogonal basis of $L^2(\R)$. In this case, the wavelet is said to be $N$-regular.
Let us introduce the notations 
\begin{equation*}
    0^d := (0,\ldots,0), \quad 1^d := (1,\ldots,1), \quad L^d = \{0,1\}^d\backslash 0^d.
\end{equation*}

An orthogonal basis of $L^2(\R^d)$ is obtained by tensorization as follows (see \cite[Chapter 2 and 3]{Meyer:1990:Ondelette_operateur} for a general construction): for $\lambda = (j,k,l)\in \Z \times \Z^d \times L^d$, let us define the tensorised wavelet 
\begin{equation*}
    \Psi_{\lambda} = \prod_{i=1}^{d} \Psi_{j,k_i}^{l_i}
\end{equation*}
with $k=(k_1,\ldots,k_D)$ and $l=(l_1,\ldots,l_D)\in L^d$, and setting $\Phi = \Psi^{0^d}$.

Thus for every $f\in L^2 (\R^d)$, one has the decomposition \eqref{eq:wlt_dec_in_L2} (observe that we use an $L^{\infty}$ normalisation for the wavelet coefficient $(c_{\lambda})_{\lambda\in\Lambda}$.)

\begin{prop}\label{prop:equi_hold_wlt_coef} \cite{Jaffard:1989:Holder_exponent_wavelet_coef,Jaffard:2004:Wavelet_techniques}
    Suppose that $\gamma>0$ and $N=\partent{\gamma+1}$. Suppose that $\Psi^{0},\Psi^{1} \in C^{N}(\R)$ is at least $N+1$-regular. Let $f:[0,1]^d \to \R$ be a locally bounded function with wavelet coefficients $\{c_{\lambda}\}_{\lambda\in\Lambda\times L^d}$, and let $x\in[0,1]^d$.

    If $f\in C^{\gamma}(x)$, then there exists a constant $M<\infty$ such that for all $\lambda=(j,k,l)\in\Lambda\times L^d$,
    \begin{equation}\label{eq:equi_hold_wlt_coef}
        |c_{\lambda}| \leq M\ 2^{-j\gamma} \big( 1+|2^j x-k| \big)^{\gamma}.
    \end{equation}
    Conversely, if \eqref{eq:equi_hold_wlt_coef} holds for $f\in \bigcup_{\varepsilon>0} C^{\varepsilon}([0,1]^d)$, then $f\in \bigcap_{\eta>0} C^{\gamma-\eta}(x)$.
\end{prop}

\subsection{Some properties of \texorpdfstring{$B^{\mu}_{p,q}(\R^d)$ and $\widetilde{B}^{\mu}_{p,q}(\R^d)$}{B mu p,q Rd}}\label{subsec:prop_inh_bes}

The spaces $\widetilde{B}^{\mu}_{p,q}(\R^d)$ and $\widetilde{B}^{\mu}_{p,q}([0,1]^d)$ are Fréchet space \cite{Barral-Seuret:2023:Besov_Space_Part_2}. A translation-invariant complete metric that induces the same topology as the family $\big\{ {\|\cdot\|_{B^{\mu^{(-\varepsilon)},\Psi}_{p,q}}} := \|\cdot\|_{L^p}+|\cdot|_{\mu^{(-\varepsilon)},p,q,\Psi} \big\}_{0<\varepsilon<\min(1,s_1)}$ can be defined by
\begin{equation}\label{eq:metric_b_mu_tilde}
    \forall f,g \in \widetilde{B}^{\mu,\Psi}_{p,q}, \quad d(f,g)=\sum_{n >\max(1,s_1^{-1})} 2^{-n} \frac{\|f-g\|_{B^{\mu^{(-1/n)},\Psi}_{p,q}}}{1+\|f-g\|_{B^{\mu^{(-1/n)},\Psi}_{p,q}}}.
\end{equation}

\begin{rem}\label{rem:incl_inh_besov} \cite[Remark 2.17]{Barral-Seuret:2023:Besov_Space_Part_2} 
    When $\mu$ satisfies $(P)$ with $0<s_1<s_2$, one has $(\mathcal{L}^d)^{\frac{s_2}{d}} \leq \mu \leq (\mathcal{L}^d)^{\frac{s_1}{d}}$ (in the sense that for every $B(x,r) \subset \R^d$, $(\mathcal{L}^d)^{\frac{s_2}{d}}(B(x,r)) \leq \mu(B(x,r)) \leq (\mathcal{L}^d)^{\frac{s_1}{d}}(B(x,r))$) inducing the embeddings
    \begin{equation}\label{eq:embed_inh_besov}
        B^{s_2+\frac{d}{p}}_{p,q}(\R^d) \hookrightarrow B^{\mu}_{p,q}(\R^d) \hookrightarrow B^{s_1+\frac{d}{p}}_{p,q}(\R^d).
    \end{equation}
    When $\alpha \geq \beta$, one has
    \begin{equation}\label{eq:embed_inh_besov_power}
        B^{\mu^{(+\alpha)}}_{p,q}(\R^d) \hookrightarrow B^{\mu^{(+\beta)}}_{p,q}(\R^d) \quad \text{and} \quad \widetilde{B}^{\mu^{(+\alpha)}}_{p,q}(\R^d) \hookrightarrow \widetilde{B}^{\mu^{(+\beta)}}_{p,q}(\R^d)
    \end{equation}
    Similarly, we can replace $\R^d$ by $[0,1]^d$ in all the previous embeddings.
\end{rem}

\section{Prevalent spectrum in \texorpdfstring{$\widetilde{B}^{\mu,\Psi}_{\infty,q}([0,1]^d)$}{b{mu,psi}{inf,q}([0,1]d)}}\label{sec:Prevalent_spectrum_inf_q}

We fix from now on $q\in[1,+\infty]$ and an almost doubling capacity $\mu\in\mathcal{C}([0,1]^d)$ with $\operatorname{supp}(\mu)=[0,1]^d$ satisfying $(P)$ and verifying the SMF (Definition \ref{defi:mf_form}).
The local behaviour of elements of $B^{\mu,\Psi}_{p,q}([0,1]^d)$ and $\widetilde{B}^{\mu,\Psi}_{p,q}([0,1]^d)$ is described via their pointwise Hölder exponent, whose study is based on the characterisation developed by Jaffard in \cite{Jaffard:2007:Wavelet_leaders}.

Consider an integer $n \geq \partent{s_2+\frac{d}{p}}+1$ and wavelets $(\Phi,\{\Psi^l\}_{l\in L^d})$ that belong to $C^n(\R^d)$ and have at least $n$ vanishing moments.

\begin{defi}\label{defi:wlt_lead}
    Let $f\in L^p_{\operatorname{loc}}(\R^d)$ for $p\in[1,+\infty]$. Let $(c_{\lambda})_{\lambda\in\Lambda \times L^d}$ be the wavelet coefficients of $f$ associated with $\Psi$. The wavelet leader of $f$ associated with $\lambda\in\Lambda$ is defined as 
    \begin{equation}
        L^{f}_{\lambda} = \sup \{ |c_{\lambda'}| : \lambda'=(j,k,l)\in \Lambda \times L^d, \lambda' \subset 3 \lambda\}.
    \end{equation}
\end{defi}

\begin{nota}
    Consider a function $f\in \widetilde{B}^{\mu,\Psi}_{p,q}([0,1]^d)$. Following Definition \ref{defi:wlt_lead}, for $x\in [0,1]^d$, the wavelet leader of $f$ at scale $j$ at $x$ is denoted by
    \begin{equation*}
        L^{f}_{j}(x) := L^{f}_{\lambda^{d}_{j}(x)},
    \end{equation*}
    with $\lambda^{d}_{j}(x) \in \Lambda^{d}_{j}$ the unique dyadic cube of generation $j$ containing $x$ (recall the definitions of Section \ref{sec:capa_prop_and_wlt_form}). The mention of $f$ in the leaders might be omitted when the function is clear from the context.
\end{nota}

\noindent The pointwise Hölder exponent can be expressed as follows \cite[Prop. 4]{Jaffard:2007:Wavelet_leaders}.

\begin{prop}\label{prop:local_dim_and_set_for_func}
    If $f \in\mathcal{C}^{\varepsilon}(R^d)$ for some $\varepsilon>0$, then for every $x_0\in\R^d$, $h_f(x_0)<r$ if and only if $\displaystyle \liminf_{j\to +\infty} \frac{\log L^{f}_{j}(x_0)}{\log (2^{-j})} <r$, and in this case 
    \begin{equation*}
        h_{f}(x_0)= \liminf_{j\to +\infty} \frac{\log L^{f}_{j}(x_0)}{\log (2^{-j})}.
    \end{equation*}
\end{prop}

\subsection{A saturating function in \texorpdfstring{$\widetilde{B}^{\mu,\Psi}_{\infty,q}$}{b mu,psi inf,q}}\label{subsec:sat_func_besov_inf_q}

To identify the prevalent properties of functions in $\widetilde{B}^{\mu,\Psi}_{p,q}([0,1]^d)$, the classical strategy is to perturb a countable family of prevalent sets with a \textit{saturating function} whose multifractal structure is the one claimed to be generic.
Let us consider the saturating function of \cite[Section 6]{Barral-Seuret:2023:Besov_Space_Part_2}.

\begin{defi}\label{defi:sat_func_besov_inf_q}
    Let us define, for every $\lambda=(j,k,l)\in\Lambda_j \times L^d$,
    \begin{equation}\label{eq:defi_sat_func:wlt_coef}
        g_{\lambda} = \frac{1}{j^{ {2}/{q}}}~\mu(\lambda^{d}) 
    \end{equation}
    with the convention $ {2}/{\infty}=0$.
    The saturating function is $\displaystyle \mathscr{G}_{q} := \sum_{\lambda \in\Lambda\times L^d} g_{\lambda} \Psi_{\lambda}$.
\end{defi}

From \cite[Lemma 6.4]{Barral-Seuret:2023:Besov_Space_Part_2}, recalling the scaling function $\tau_{\mu}$ \eqref{eq:scaling_function_tau_mu}, one has:
\begin{prop}
    We have $\mathscr{G}_{q} \in \widetilde{B}^{\mu,\Psi}_{\infty,q}([0,1]^d)$ and
    $\sigma_{\mathscr{G}_{q}} = \tau^{*}_{\mu} = \sigma_{\mu}$.
\end{prop}

One actually has that $\mathscr{G}_{q}\in B^{\mu,\Psi}_{\infty,q}([0,1]^d) \subset \widetilde{B}^{\mu,\Psi}_{\infty,q}([0,1]^d)$ directly from \eqref{eq:snorm_inh_besov} which gives that $|\mathscr{G}_{q}|_{\mu,\infty,q,\Psi}^q=\sum_{j=1}^{\infty} j^{-2}$ that converges. The singularity spectrum is deduced through the key property that for every $x\in(0,1)^d$, $h_{\mathscr{G}_{q}}(x)=h_{\mu}(x)$.

\subsection{Lower bound of the pointwise Hölder exponent in \texorpdfstring{$\widetilde{B}^{\mu,\Psi}_{\infty,q}([0,1]^d)$}{b{mu,psi}{inf,q}([0,1]d)}}\label{subsec:low_bound_hold_exp_b_tilde}

Let us begin with the following result.
\begin{lemme}\label{lemma:min_exp_besov_inf_q}
    For all $f\in \widetilde{B}^{\mu,\Psi}_{\infty,q}([0,1]^d)$, for every $x\in[0,1]^{d}$, $h_{f}(x) \geq \underline{h}_{\mu}(x)$.
\end{lemme}
\begin{proof}
    Let us write the decomposition $f= c_{0,0^d,0^d} \Psi_{0,0^d,0^d} + \sum_{\lambda \in \Lambda\times L^d} c^f_{\lambda} \Psi_{\lambda}$.
    Let $\varepsilon\in(0,\min(1,s_1))$.
    Recall Definition \ref{defi:inh_besov_wlt}. Since $f\in B^{\mu^{(-\varepsilon)},\Psi}_{\infty,q}([0,1]^d)$, for every $j\in\N$, $\varepsilon^{\mu^{(-\varepsilon)},\infty}_j < + \infty$ and one has
    \begin{equation*}
        \sup_{j\in\N}\sup_{\lambda \in \Lambda_{j}\times L^d} \bigg|\frac{c^f_{\lambda}}{\mu(\lambda)~2^{j\varepsilon}}\bigg| \leq C.
    \end{equation*}
    For a fixed $x \in [0,1]^d$, one has
    \begin{equation*}
        L^{f}_{j}(x)
        = \sup_{\lambda \subset 3 \lambda^{d}_{j}(x)} \big|c^{f}_{\lambda}\big| 
        \leq C\sup_{\lambda_{j'} \subset 3 \lambda^{d}_{j}(x)} |\mu(\lambda_{j'})| 2^{j' \varepsilon}.
    \end{equation*}
    Recall that $\mu$ satisfy property $(P)$ with $0<s_1<s_2$. Hence, setting $M = \partent{\frac{s_2}{s_1}}+1$, \eqref{eq:P1} gives that for every $j'>Mj$, if $\lambda_{j'} \subset 3 \lambda^{d}_{j}(x)$, then $\mu(\lambda_{j'})2^{j' \varepsilon} \leq \widetilde{C}\,2^{-j' (s_1-\varepsilon)} \leq \widetilde{C}\,2^{-j (s_2-\varepsilon)} \leq \widetilde{C}^2\,\mu(\lambda^{d}_{j}(x))2^{j \varepsilon}$.
    For $j\leq j' \leq Mj$, using that $\mu$ is a capacity, one has $\mu(\lambda_{j'})2^{j' \varepsilon} \!\leq\! \mu(\lambda_{j'})2^{Mj \varepsilon} \!\leq\! \mu(3\lambda^{d}_{j}(x))2^{Mj \varepsilon}$.
    The almost-doubling property of $\mu$ gives that there exists a non-increasing $\phi_{\mu}:[0,1) \to \R_+$, such that $\lim_{r\to 0^+} \frac{\phi_{\mu}(r)}{\log(r)}=0$, for which \eqref{eq:alm_doub_rad} applies. Denote $y$ the center of $\lambda^{d}_{j}(x)$ so $\lambda^{d}_{j}(x) = B\big(y,2^{-(j+1)}\big)$ and remark that $3\lambda^{d}_{j}(x) \subset B\big(y,2^{-(j-1)}\big)$. Considering $r=2^{-j}$, one has $\mu(3\lambda^{d}_{j}(x)) \leq \mu(B(y,2^{-(j-1)})) \leq e^{\phi_{\mu}(2^{-(j+1)})+\phi_{\mu}(2^{-j})} \mu(\lambda^{d}_{j}(x)) \leq e^{2\phi_{\mu}(2^{-(j+1)})} \mu(\lambda^{d}_{j}(x))$.
    Combining this with the last two inequalities gives, for large $j$s,
    \begin{align*}
        L^{f}_{j}(x) 
        &\leq C\max(\widetilde{C}^22^{j\varepsilon}, e^{2\phi_{\mu}(2^{-(j+1)})}2^{Mj \varepsilon})\, \mu(\lambda^{d}_{j}(x))\\
        &\leq C\, e^{2\phi_{\mu}(2^{-(j+1)})}2^{Mj \varepsilon} \mu(\lambda^{d}_{j}(x)).
    \end{align*}
    Hence, for all $j\in\N$ large enough
    \begin{align*}
        \log L^{f}_{j}(x) &\leq \log C + 2\phi_{\mu}\big(2^{-(j+1)}\big) + \log\mu(\lambda^{d}_{j}(x)) + \log\big(2^{Mj \varepsilon}\big)\\
        \text{and}~ \frac{\log L^{f}_{j}(x)}{\log (2^{-j})} &\geq \frac{\log C}{\log (2^{-j})} + 2\frac{\phi_{\mu}\big(2^{-(j+1)}\big)}{\log (2^{-j})} + \frac{\log \mu(\lambda^{d}_{j}(x))}{\log (2^{-j})} - M\varepsilon.
    \end{align*}
    Taking the $\liminf_{j\in\N}$ on both sides gives $\displaystyle h_{f}(x) \geq \underline{h}_{\mu}(x) - M\varepsilon$.
    Since this holds for all $\varepsilon\in(0,\min(1,s_1))$, one concludes that $\displaystyle h_{f}(x) \geq \underline{h}_{\mu}(x)$.
\end{proof}

\begin{rem}
    In Lemma \ref{lemma:min_exp_besov_inf_q}, one can replace $\widetilde{B}^{\mu,\Psi}_{\infty,q}([0,1]^d)$ by $B^{\mu,\Psi}_{\infty,q}([0,1]^d)$, by simply setting $\varepsilon$ to zero in all the computations.
\end{rem}

\subsection{Proof of Theorem \ref{theo:prev_spect_inh_besov_inf_q}}

The key result to obtain the prevalent value of the singularity spectrum of functions in $\widetilde{B}^{\mu,\Psi}_{\infty,q}([0,1]^d)$ in Theorem \ref{theo:prev_spect_inh_besov_inf_q} is the following:

\begin{prop}\label{prop:prev_set_besov_inf_q}
    The set 
    $$\mathcal{F}^{\mu} = \big\{ f\in \widetilde{B}^{\mu,\Psi}_{\infty,q}([0,1]^d):\forall x \in [0,1]^d,~h_{f}(x) \leq \overline{h}_{\mu}(x)\big\}$$
    is prevalent in $\widetilde{B}^{\mu,\Psi}_{\infty,q}([0,1]^d)$.
\end{prop}

The proof of Proposition \ref{prop:prev_set_besov_inf_q} is postponed to Subsection \ref{subsubsec:proof_prop_prev_set}. For the moment, we admit it and explain how this allows us to deduce Theorem \ref{theo:prev_spect_inh_besov_inf_q}. 

\begin{proof}[Proof of Theorem \ref{theo:prev_spect_inh_besov_inf_q}]
    Using Lemma \ref{lemma:min_exp_besov_inf_q} and Proposition \ref{prop:prev_set_besov_inf_q}, for all $f\in\mathcal{F}^{\mu}$, for every $x\in[0,1]^{d}$, 
    \begin{equation}\label{eq:bound_h_f}
        \underline{h}_{\mu}(x) \leq h_{f}(x) \leq \overline{h}_{\mu}(x).
    \end{equation}
    
    Fix $h \in\R_{+}$. 
    
    \underline{Step 1:} we start with \textit{the upper bound} $\sigma_{f} \leq \sigma_{\mu}$. By \eqref{eq:bound_h_f}, one gets
    \begin{align*}
        E_{f}(h) &\subset \{x \in \R^{d}:h_{f}(x) \leq h\} \subset \{x \in \R^{d}:\underline{h}_{\mu}(x) \leq h\} = \underline{E}_{\mu}^{\leq}(h), \\
        E_{f}(h) &\subset \{x \in \R^{d}:h_{f}(x) \geq h\} \subset \{x \in \R^{d}:\overline{h}_{\mu}(x) \geq h\} = \overline{E}_{\mu}^{\geq}(h).
    \end{align*}
    We now apply Proposition \ref{prop:prop_capa_SMF} to $\mu$ which satisfies the SMF, therefore
    \begin{itemize}[noitemsep]
        \item if $h \leq \tau'_{\mu}(0^{+})$, \ref{item:prop_capa_SMF:dim_E_leq} gives 
        \begin{equation*}
            \sigma_{f}(h) = \dimH(E_{f}(h)) \leq \dimH(\underline{E}_{\mu}^{\leq}(h)) = \tau_{\mu}^{*}(h) = \sigma_{\mu}(h),
        \end{equation*}
        \item if $ \tau'_{\mu}(0^{+}) \leq h \leq \tau'_{\mu}(0^{-})$, with \ref{item:prop_capa_SMF:sigma_is_d}, one has $\displaystyle \sigma_{f}(h) \leq d = \sigma_{\mu}(h)$,
        \item if $h \geq \tau'_{\mu}(0^{-})$, \ref{item:prop_capa_SMF:dim_E_geq} gives 
        \begin{equation*}
            \sigma_{f}(h) = \dimH(E_{f}(h)) \leq \dimH(\overline{E}_{\mu}^{\geq}(h)) = \tau_{\mu}^{*}(h) = \sigma_{\mu}(h).
        \end{equation*}
    \end{itemize}
    Hence, one gets the upper bound.

    \underline{Step 2:} let us move to \textit{the lower bound} $\sigma_{f} \geq \sigma_{\mu}$. 
    For all $x\in E_{\mu}(h)$, by Definition \ref{defi:local_dim_and_set}, one has $\displaystyle \underline{h}_{\mu}(x) = \overline{h}_{\mu}(x) = h_{\mu}(x)$. 
    Thus, by \eqref{eq:bound_h_f}, one deduces $h_{\mu}(x) = h_{f}(x)$ and $E_{\mu}(h) \subset E_{f}(h)$. Finally, Proposition \ref{prop:prop_capa_SMF} gives $\sigma_{f}(h) = \dimH(E_{f}(h)) \geq \dimH(E_{\mu}(h)) = \sigma_{\mu}(h)$.
\end{proof}

\subsubsection{Prevalence property of the auxiliary set \texorpdfstring{$\mathcal{F}^{\mu}$}{F of xi}}\label{subsubsec:proof_prop_prev_set}

We want to show that $(\mathcal{F}^{\mu})^\complement = \widetilde{B}^{\mu,\Psi}_{\infty,q}([0,1]^d) \backslash \mathcal{F}^{\mu}$ is shy.
Observe that
 $$(\mathcal{F}^{\mu})^{\complement} = \big\{ f\in \widetilde{B}^{\mu,\Psi}_{\infty,q}([0,1]^d):\exists x \in [0,1]^d,~ h_{f}(x) > \overline{h}_{\mu}(x) \big\}.$$

We simplify the problem by including $(\mathcal{F}^{\mu})^{\complement}$ in a countable union of simpler ancillary sets.
For $p\in\N^{*}$, $\gamma>0$ and $M=(M_1,M_2)\in(\N^{*})^2$, let
\begin{equation*}
    \mathcal{O}^{p,\gamma,M} = \left\{f\in \widetilde{B}^{\mu,\Psi}_{\infty,q}([0,1]^d):
    \begin{array}{c}
        \exists x \in [0,1]^d,\forall \lambda_j = (j,k,l)\in\Lambda\times L^d,\\
        \big|c^{f}_{\lambda_j}\big| \leq M_1\ 2^{-j\gamma} \big(1+|2^j x-k|\big)^{\gamma}\\
        \forall j>M_2,~\mu(\lambda_j(x)) \geq 2^{-j(\gamma-1/p)}
    \end{array}\right\}.
\end{equation*}
Remembering Proposition \ref{prop:equi_hold_wlt_coef}, the sets $\mathcal{O}^{p,\gamma,M}$ shall be understood as the set of functions in $\widetilde{B}^{\mu,\Psi}_{\infty,q}([0,1]^d)$ such that for some $x\in [0,1]^d$, $f \in C^{\gamma}(x)$ and simultaneously the $\mu$-mass of the dyadic cubes $\lambda_j(x)$ is bounded below. This two properties will allow us to control the ratio $\frac{c^{f}_{\lambda_j(x)}}{\mu(\lambda_j(x))}$.

\begin{prop}\label{prop:incl_cFxi_in_O}
    One has $(\mathcal{F}^{\mu})^\complement \subset \bigcup_{p\in\N^*} \bigcup_{\gamma\in\Q_{+}} \bigcup_{M\in(\N^*)^2} \mathcal{O}^{p,\gamma,M}$.
\end{prop}
\begin{proof}
    We write $\displaystyle (\mathcal{F}^{\mu} )^{\complement} = \bigcup_{p\in\N^*} \mathcal{F}_{p}$, where $\mathcal{F}_{p} = \big\{ f\in \widetilde{B}^{\mu,\Psi}_{\infty,q}([0,1]^d):\exists x \in [0,1]^d,~h_{f}(x) - \overline{h}_{\mu}(x) > \frac{2}{p} \big\}$.
    Supposing that $f\in \mathcal{F}_{p}^{\mu}$, there exist $\gamma \in \Q_{+}$ and $x \in [0,1]^d$ such that $h_{f}(x) > \gamma$ and $\overline{h}_{\mu}(x) \leq \gamma- \frac{2}{p}$. Using \eqref{eq:equi_hold_wlt_coef} of Proposition \ref{prop:equi_hold_wlt_coef}, one has
    \begin{align*}
        h_{f}&(x) > \gamma 
        \Longrightarrow f \in C^{\gamma}(x) \\
        &\Longrightarrow \exists M_{1}\in\N,\forall \lambda_j = (j,k,l)\in\Lambda\times L^d, \big|c^{f}_{\lambda_j}\big| \leq M_{1}\ 2^{-j\gamma} \big(1+|2^j x-k|\big)^{\gamma}.
    \end{align*}
    In addition, by \eqref{eq:h_mu_dyadic}, one has
    \begin{align*}
        \overline{h}_{\mu}(x) \leq \gamma - \frac{2}{p}
        &\Longrightarrow \exists M_2\in\N,~\sup_{j>M_2} \frac{\log \mu(\lambda_j(x))}{\log 2^{-j}} \leq \gamma-\frac{1}{p} \\
        &\Longrightarrow \exists M_2\in\N,~\forall j>M_2,~\mu(\lambda_j(x)) \geq 2^{-j(\gamma-1/p)}.
    \end{align*}
    Hence, $f\in\mathcal{O}^{p,\gamma,M}$ with $M=(M_1,M_2)$. 
    One gets that if $f\notin \mathcal{F}^{\mu}$, then there exists $p\in\N^*$ such that $f \in\mathcal{F}_{p}^{\mu}$. Then, for $p\in\N^*$, there are $\gamma\in\Q_{+}$ and an $M=(M_1,M_2)\in(\N^*)^2$ such that $f\in \mathcal{O}^{p,\gamma,M}$.
\end{proof}

To use the prevalence theory, it remains to prove the universal measurability of $\mathcal{O}^{p,\gamma,M}$. Analytic sets in Polish spaces are defined as continuous images of Borel sets. Unfortunately, as we work with $\widetilde{B}^{\mu,\Psi}_{\infty,q}([0,1]^d)$ which is not separable, we cannot use the Polish spaces characterisation and have to consider a more general definition, from Choquet \cite{Choquet:1954:theory_of_capacities}, valid for any Hausdorff topological space $X$ endowed with a Borel $\sigma$-algebra $\mathcal{B}(X)$. 
Choquet's capacitability theorem \cite[Theorem 30.1]{Choquet:1954:theory_of_capacities} with the following definition of analyticity \cite[Theorem 5.1]{Choquet:1954:theory_of_capacities} give the right setting, as was observed in \cite{Aubry-Maman-Seuret:2013:Traces_besov_results}.
\begin{defi}\label{defi:ana_haus_space}
    For a compact topological space $K$, the collection of its closed subsets is written $\mathcal{K}$. Let $(\mathcal{B}(X) \times \mathcal{K})_{\sigma \delta}$ be the collection of countable intersections of countable unions of sets that are Cartesian products of a Borel set in $X$ and a closed set in $K$.
    
    A set $A\subset X$ is said to be analytic if there exist a compact space $K$ and $\mathcal{T}\in(\mathcal{B}(X)\times\mathcal{K})_{\sigma \delta}$ such that $A= \pi(\mathcal{T})$
    with $\pi:X\times K \to X$ the canonical projection map $\pi(x,y)=x$.
\end{defi}
Then Choquet's capacitability theorem yields the following theorem.
\begin{theo}\label{theo:ana_to_uni_meas}
    Every analytic set of a Hausdorff topological space is universally measurable.
\end{theo}

Theorem \ref{theo:ana_to_uni_meas} states the equivalence between the universal measurability of $\mathcal{O}^{p,\gamma,M}$ and its analyticity.

\begin{prop}\label{prop:uni_meas_O}
    Let $p\in\N^*$, $M\in(\N^{*})^2$ and $\gamma\in\Q_{+}$. The set $\mathcal{O}^{p,\gamma,M}$ is universally measurable in $\widetilde{B}^{\mu,\Psi}_{\infty,q}([0,1]^d)$.
\end{prop}
\begin{proof}
    We are going to prove the analyticity of $\mathcal{O}^{p,\gamma,M}$ by taking the right sets $X$, $K$ and $\mathcal{T}$ in Definition \ref{defi:ana_haus_space}, then writing it as $\pi(\mathcal{T})$.
    Denote $X:= \big(\widetilde{B}^{\mu,\Psi}_{\infty,q}([0,1]^d),d \big)$ with $d$ defined in \eqref{eq:metric_b_mu_tilde} and $K:=([0,1]^{d},\|\cdot\|_{\infty})$.
    \medskip

    \underline{Step 1:} For $\lambda_j=(j,k,l)\in\Lambda \times L^d$, let us introduce the maps $\phi_{\lambda_j}, \phi : X\times K \to \R$ defined as
    \begin{align*}
        \phi_{\lambda_j} : (f,x) &\mapsto M_1~2^{-j\gamma} \big(1+|2^j x-k|\big)^{\gamma} - \big|c^{f}_{\lambda_j}\big| \\
        \phi : (f,x) &\mapsto \inf_{\lambda\in\Lambda\times L^d} \phi_{\lambda} (f,x).
    \end{align*}
    \begin{lemme}
        Each $\phi_{\lambda_j}$ is continuous on $X\times K$. 
    \end{lemme}
    \begin{proof}
        Let $\varepsilon>0$. 
        Recalling \eqref{eq:metric_b_mu_tilde}, for any $f,g \in X$ such that $d(f,g) < \varepsilon$, one gets that, for any $n:= \partent{\max(1,s_1^{-1})}+1$, 
        \begin{equation*}
            2^{-n} \frac{\|f-g\|_{B^{\mu^{(-1/n)},\Psi}_{\infty,q}}}{1+\|f-g\|_{B^{\mu^{(-1/n)},\Psi}_{\infty,q}}} < \varepsilon.
        \end{equation*}
        This induces that $\|f-g\|_{B^{\mu^{(-1/n)},\Psi}_{\infty,q}} < 2^{n} \varepsilon$. By definition of $\|\cdot\|_{B^{\mu^{(-\varepsilon)},\Psi}_{\infty,q}}$ in Section \ref{subsec:prop_inh_bes}, $\|f-g\|_{L^\infty} < 2^{n} \varepsilon$.
        
        The application $f \in X \mapsto c^{f}_{\lambda_j}=\pscal{f}{\Psi_{\lambda_j}} \in \R$ is a scalar product between $f$ and a fixed wavelet $\Psi_{\lambda_j}$. With the classical remark that $\pscal{f}{\Psi_{\lambda_j}} \leq \|f\|_{L^\infty} \|\Psi_{\lambda_j}\|_{L^1}$, $f \mapsto c^{f}_{\lambda_j}$ is continuous. 
        So the mapping $(f,x) \mapsto c^{f}_{\lambda_j}$ is continuous and so is $\phi_{\lambda_j}$ as a sum of continuous functions.
    \end{proof}

    \underline{Step 2:} Denote
    \begin{equation}\label{eq:def_D_mu}
        \mathcal{D}_{\mu}(j,\gamma-1/p) := \bigg\{ \lambda_j\in\Lambda_{j} : \frac{\log\mu(\lambda_j)}{\log(2^{-j})} \leq \gamma - \frac{1}{p}\bigg\}.
    \end{equation}
    Let us introduce $\mathcal{T} := \phi^{-1} ([0,+\infty)) \bigcap \Big(X \times \bigcap_{j\geq M_2} \bigcup_{\lambda\in\mathcal{D}_{\mu}(j,\gamma-1/p)} \lambda \Big)$.
    \begin{lemme}\label{lemma:T_belong_sigma_delta}
        The set $\mathcal{T}$ belongs to $(\mathcal{B}(X) \times \mathcal{K})_{\sigma \delta}$.
    \end{lemme}
    \begin{proof}
        First, it is straightforward from the definition that, for all $\lambda\in\Lambda$, $X \times \lambda \in \mathcal{B}(X) \times \mathcal{K}$, so taking a countable intersection of compact sets,
        \begin{equation*}
            X \times \bigcap_{j\geq M_2} \bigcup_{\lambda\in\mathcal{D}_{\mu}(j,\gamma-1/p)} \lambda \in (\mathcal{B}(X) \times \mathcal{K})_{\sigma \delta}.
        \end{equation*}
        
        Let $n\in\N$ and $m\in \Z$ and consider
        \begin{equation*}
            F_{\lambda_j}(n,m) := \{f\in \widetilde{B}^{\mu,\Psi}_{\infty,q}([0,1]^d):c^f_{\lambda_j} \in 2^{-n} [m, m+1]\}.
        \end{equation*}
        Since $2^{-n} [m, m+1]$ is a closed interval, $F_{\lambda_j}(n,m)$ is closed as the preimage of a closed set by the continuous function $f\mapsto c_{\lambda_j}^{f}$.
        Next, let us define
        \begin{equation*}
            X_{\lambda_j}(n,m)
            := \bigg\{ x\in[0,1]^d:\sup_{f\in F_{\lambda_j}(n,m)} \phi_{\lambda_j}(f,x) \geq 0 \bigg\}.
        \end{equation*}
        Observe that when $f$ ranges in $\widetilde{B}^{\mu,\Psi}_{\infty,q}([0,1]^d)$, all possible values of $c^f_{\lambda_j} \in 2^{-n} [m, m+1]$ are reached. Hence, one can rewrite $X_{\lambda_j}(n,m)$ as
        \begin{equation*}
            X_{\lambda_j}(n,m) 
            = \Big\{ x\in[0,1]^d : \max_{c \in 2^{-n} [m, m+1]} M_1 \,2^{-j\gamma} \big(1+|2^j x-k|\big)^{\gamma} - |c| \geq 0 \Big\}.
        \end{equation*}
        The $\max$ can be in turn rewritten 
        \begin{equation*}
            X_{\lambda_j}(n,m) 
            = \bigg\{ x\in[0,1]^d :
            \begin{array}{c}
            \exists c\in 2^{-n} [m, m+1], \\
            M_1 2^{-j\gamma} \big(1+|2^j x-k|\big)^{\gamma} - |c| \geq 0
            \end{array} \bigg\}.
        \end{equation*}
        Denoting the mapping $v: (x,c) \in\R^d\times\R \mapsto M_1\,2^{-j\gamma} \big(1+|2^j x-k|\big)^{\gamma} - |c| \in \R$, one has, for $\widetilde{\pi}(x,c)=x$ the projection on the first coordinate,
        \begin{equation*}
            X_{\lambda_j}(n,m) 
            = \widetilde{\pi} \big( v^{-1}(\R_{+}) \cap ([0,1]^d\times 2^{-n} [m, m+1]) \big).
        \end{equation*}
        The mapping $v:\R^d\times\R \to \R$ is continuous, so $v^{-1}(\R_{+}) \cap ([0,1]^d\times 2^{-n} [m, m+1])$ is compact. As the projection is continuous, $X_{\lambda_j}(n,m)$ is also compact.
        \begin{lemme}
            One has $\phi_{\lambda_j}^{-1} ([0,+\infty)) = \bigcap_{n\in\N} \bigcup_{m\in \Z} F_{\lambda_j}(n,m) \times X_{\lambda_j}(n,m)$.
        \end{lemme}
        \begin{proof}
            We first prove that $\phi_{\lambda_j}^{-1} ([0,+\infty)) \!\subset\! \bigcap_{n\in\N} \bigcup_{m\in \Z} \!F_{\lambda_j}\!(n,m) \times X_{\lambda_j}\!(n,m)$.

            Let $(f,x)\in \phi_{\lambda_j}^{-1} ([0,+\infty))$. Then for every $n\in\N$, there is an integer $m\in\Z$ such that $c^f_{\lambda_j} \in 2^{-n} [m, m+1]$ and so $f\in F_{\lambda_j}(n,m)$. In addition, $\phi_{\lambda_j}(f,x) \geq 0$, so $\sup_{g\in F_{\lambda_j}(n,m)} \phi_{\lambda_j}(g,x) \geq 0$ and $x\in X_{\lambda_j}(n,m)$.
            Hence $(f,x) \in F_{\lambda_j}(n,m) \times X_{\lambda_j}(n,m)$.
            \medskip

            Next, let us prove $\phi_{\lambda_j}^{-1} ([0,+\infty)) \supset \bigcap_{n\in\N} \bigcup_{m\in \Z} F_{\lambda_j}(n,m) \times X_{\lambda_j}(n,m)$.

            Let $(f,x)\in \bigcap_{n\in\N} \bigcup_{m\in \Z} F_{\lambda_j}(n,m) \times X_{\lambda_j}(n,m)$. Then for all $n\in\N$, there is an $m\in\Z$ such that $f\in F_{\lambda_j}(n,m)$ and $x\in X_{\lambda_j}(n,m)$.
            
            One has $c^f_{\lambda_j} \in 2^{-n} [m, m+1]$ and $\sup_{g\in F_{\lambda_j}(n,m)} \phi_{\lambda_j}(g,x) \geq 0$. Hence
            \begin{align*}
                \bigg| \sup_{g\in F_{\lambda_j}(n,m)} \phi_{\lambda_j}(g,x) - \phi_{\lambda_j}(f,x) \bigg|
                &\leq \sup_{g\in F_{\lambda_j}(n,m)} |\phi_{\lambda_j}(g,x) - \phi_{\lambda_j}(f,x)\ \\
                &\leq \sup_{g\in F_{\lambda_j}(n,m)} \big|c^{g}_{\lambda_j} - c^{f}_{\lambda_j}\big| \leq 2^{-n}.
            \end{align*}
            As, this result holds for all $n\in\N$, one deduces that $\phi_{\lambda_j}(f,x) \geq 0$
        \end{proof}
        Hence, $\phi_{\lambda_j}^{-1} ([0,+\infty)) \in (\mathcal{B}(X) \times \mathcal{K})_{\sigma \delta}$. We deduce that $\phi^{-1} ([0,+\infty)) = \bigcap_{\lambda\in \Lambda\times L^d} \phi_{\lambda}^{-1} ([0,+\infty))$ also belongs to $(\mathcal{B}(X) \times \mathcal{K})_{\sigma \delta}$. 
        This proves that the set $\mathcal{T}$ indeed belongs to $(\mathcal{B}(X) \times \mathcal{K})_{\sigma \delta}$.
    \end{proof}

    \underline{Step 3:} By Definition \ref{defi:ana_haus_space}, the canonical projection map $\pi:(x,y)\in X\times K \to x\in X$ gives that 
    \begin{equation*}
        \pi(\mathcal{T}) := \bigg\{ f\in X : 
        \begin{array}{c}
            \exists x\in K,(f,x)\in\phi^{-1} ([0,+\infty)),\\
            \forall j>M_2,~\mu(\lambda_j(x)) \geq 2^{-j(\gamma-1/p)}
        \end{array}\bigg\}
    \end{equation*}
    is analytic in the space $(X,\mathcal{B}(X))$.
    \medskip

    \underline{Step 4:} To conclude, notice that one can rewrite $\mathcal{O}^{p,\gamma,M}=\pi(\mathcal{T})$.
    This yields that $\mathcal{O}^{p,\gamma,M}$ is analytic and thus universally measurable.
\end{proof}

As the sets $\mathcal{O}^{p,\gamma,M}$ are universally measurable, prevalence theory is now used to show the shyness of $(\mathcal{F}^{\mu})^\complement$ through the shyness of $\mathcal{O}^{p,\gamma,M}$. To do so, the classical method is to use a finite-dimensional subspace $\mathcal{P}$ and the Lebesgue measure restricted to a compact subset of $\mathcal{P}$ as a transverse measure.


Let $d_1$ be an integer such that 
\begin{equation}\label{eq:def_d_1}
    d < d_1 / 2p.
\end{equation}

\begin{defi}\label{defi:sat_func_multi_dim_subspace}
    Consider the saturating function $\mathscr{G}_{q}$ whose wavelet coefficients are $g_{\lambda}$ (see \eqref{eq:defi_sat_func:wlt_coef}). 
    Let $\mathcal{P} = \operatorname{Vect}_{i=1,\ldots,d_1}(\mathscr{G}_{q}^{(i)}) \subset \widetilde{B}^{\mu,\Psi}_{\infty,q}([0,1]^d)$ where $\mathscr{G}_{q}^{(i)}$ is the wavelet series whose wavelet coefficients $g^{(i)}_{\lambda}$ are defined in the following way: for every $\lambda=(j,k,l)\in \Lambda\times \{0,1\}^d$,
    \begin{equation}\label{eq:defi:sat_func_multi_dim_subspace}
        g^{(i)}_{\lambda} = 
        \begin{cases}
            \displaystyle g_{\lambda} &\text{if } j\pmod{d_1} = i,\\
            0 &\text{else}.
        \end{cases}
    \end{equation}
\end{defi}
Take an arbitrary $f\in \widetilde{B}^{\mu,\Psi}_{\infty,q}([0,1]^d)$, define the function $f^{\beta}$ by
\begin{equation}
    \begin{array}{ccc}
        \R^{d_1} & \to & f+\mathcal{P} \\
        \beta & \mapsto & f^{\beta} = f+ \sum_{i=1}^{d_1} \beta_i \cdot \mathscr{G}_{q}^{(i)}.
    \end{array}
\end{equation}

Denote $c^{f}_{\lambda}$ (resp. $c_{\lambda}^{\beta}$) the wavelet coefficient of $f$ (resp. $f^{\beta}$) and $L^{f}_{\lambda}$ (resp. $L_{\lambda}^{\beta}$) the associated wavelet leader.

\begin{prop}\label{prop:beta_B_meas_0}
    For any function $f\in \widetilde{B}^{\mu,\Psi}_{\infty,q}([0,1]^d)$, the set $\{ \beta\in\R^{d_1} : f^{\beta} \in \mathcal{O}^{p,\gamma,M} \}$ has $d_1$-dimensional Lebesgue measure $\mathcal{L}_{d_1}$ zero.
\end{prop}
\begin{proof}
    Let us rewrite $\displaystyle \big\{ \beta\in\R^{d_1} : f^{\beta} \in \mathcal{O}^{p,\gamma,M} \big\}$ as
    \begin{equation*}
        \mathcal{B}_f := 
        \left\{\beta\in\R^{d_1}:
        \begin{array}{c}
            \exists x \in [0,1]^d,\forall \lambda_j = (j,k,l)\in\Lambda\times L^d,\\
            \big|c^{\beta}_{\lambda_j}\big| \leq M_1~2^{-j\gamma} \big(1+|2^j x-k|\big)^{\gamma},\\
            \forall j>M_2,~\mu(\lambda_j(x)) \geq 2^{-j(\gamma-1/p)}
        \end{array}\right\}.
    \end{equation*}
    We start with a technical lemma.
    \begin{lemme}\label{lemma:B_map_meas}
        The mapping $\displaystyle \left\{ \begin{array}{ccc}
        \R^{d_1} & \longrightarrow & \R \\
        \beta & \longmapsto & \mathds{1}_{\mathcal{B}}(\beta)
        \end{array} \right.$ is measurable.
    \end{lemme}
    \begin{proof}
        We will use the notation and function from the proof of Lemma \ref{lemma:T_belong_sigma_delta}.
        
        Let $\displaystyle \Phi : (\beta,x) \mapsto \inf_{\lambda \in\Lambda\times L^d} \phi_{\lambda}(f^{\beta},x)$. The map $\Phi$ is Borel on $\R^{d_1} \times [0,1]^{d}$ as an infimum of countably many continuous functions. One sees that 
        \begin{equation*}
            \mathcal{B}_f= \bigg\{ \beta \in \R^{d_1}:
            \begin{array}{c}
                \exists x\in [0,1]^d,~\Phi(\beta,x)\geq 0, \\
                \forall j>M_2,~\mu(\lambda_j(x)) \geq 2^{-j(\gamma-1/p)}
            \end{array}\bigg\}.
        \end{equation*}
        This set can be rewritten as
        \begin{equation}\label{eq:meas_proj_borel_set}
            \mathcal{B}_f= \widetilde{\pi}\bigg( \Phi^{-1}([0,+\infty)) \bigcap \bigg( \R^{d_1} \times \bigcap_{j\geq M_2} \bigcup_{\lambda\in\mathcal{D}_{\mu}(j,\gamma-1/p)} \lambda \bigg)\bigg),
        \end{equation}
        where $\widetilde{\pi}(\beta,x)=\beta$ is the canonical projection on the first coordinate and $\mathcal{D}_{\mu}(j,\gamma-1/p)$ is defined by \eqref{eq:def_D_mu}. Since the set in brackets in \eqref{eq:meas_proj_borel_set} is a Borel set, $\mathcal{B}_f$ is analytic by Definition \ref{defi:ana_haus_space}. By Theorem \ref{theo:ana_to_uni_meas}, it is also Lebesgue-measurable.
    \end{proof}
    
    We are ready to prove that the set $\mathcal{B}_f$ has $d_1$-dimensional Lebesgue measure $\mathcal{L}_{d_1}$ zero.
    For any $\lambda_0 := (j_0,k_0,l_0)\in\Lambda$, denote
    \begin{equation*}
        \mathcal{B}_{f,\lambda_0} = \left\{ \beta\in\R^{d_1} :
        \begin{array}{c}
            \exists x \in \lambda_0,\forall \lambda_j = (j,k,l)\in\Lambda\times L^d, \\
            \big|c^{\beta}_{\lambda_j}\big| \leq M_1~2^{-j\gamma} \big(1+|2^j x-k|\big)^{\gamma},\\
            \forall j>M_2,~\mu(\lambda_j(x)) \geq 2^{-j(\gamma-1/p)}
        \end{array}\right\}.
    \end{equation*}
    so that
    \begin{equation}\label{eq:dec_B}
        \mathcal{B}_f = \limsup_{j_0\to \infty} \bigcup_{(j_0,k_0,l_0)\in\Lambda_{j_0}^d} \mathcal{B}_{f,\lambda_0} = \bigcap_{J\in\N^*} \bigcup_{j_0\geq J} \bigcup_{(j_0,k_0,l_0)\in\Lambda_{j_0}^d} \mathcal{B}_{f,\lambda_0}.
    \end{equation}
    Let us show that $\mathcal{L}_{d_1}(\mathcal{B}_{f}) = 0$ by establishing an upper bound for each $\mathcal{L}_{d_1}(\mathcal{B}_{f,\lambda_0})$ and using the Borel-Cantelli lemma.
    One can ensure that $j_0$ is large enough so that $j_0 \geq M_2 + d_1$.
    Suppose that $\beta$ and $\widetilde{\beta}$ belong to a same $\mathcal{B}_{f,\lambda_0}$. 
    There exist $x_{\beta}, x_{\widetilde{\beta}}\in \lambda_0 \subset \R^d$ such that for every $j\in\N^*$, 
    \begin{equation}\label{eq:maj_c_j_beta}
        \big|c^{\beta}_{\lambda_{j}(x_{\beta})}\big| \leq M_1\,2^{-j\gamma} \big(1+|2^j x_{\beta}-k_{j,\beta}|\big)^{\gamma} 
    \end{equation} 
    and similarly for $c^{\widetilde{\beta}}_{\lambda_{j}(x_{\widetilde{\beta}})}$ with $k_{j,\beta}$ and $k_{j,\widetilde{\beta}}$ the multi-integers associated respectively to $\lambda_{j}(x_{\beta})$ and $\lambda_{j}(x_{\widetilde{\beta}})$.

    For every $j\in\{j_0-d_1,\ldots,j_0-1\}$, one has $\lambda_{j}(x_{\beta}) = \lambda_{j}(x_{\widetilde{\beta}})$. Let us denote $\widetilde{\lambda}_j=\lambda_{j}(x_{\beta})$. 
    
    Using $x_{\beta}, x_{\widetilde{\beta}}\in \lambda_0 \subset \widetilde{\lambda}_j$ and \eqref{eq:maj_c_j_beta}, one deduces that there is a constant $C>0$ depending on $M, j_0, d_1$ and $\gamma$ such that, for every $j\in\{j_0-d_1,\ldots,j_0-1\}$, $\big|c^{\beta}_{\widetilde{\lambda}_j}\big|\leq C\, 2^{-j\gamma}$ and $\big|c^{\widetilde{\beta}}_{\widetilde{\lambda}_j}\big| \leq C\, 2^{-j\gamma}$.
    
    From all this, we deduce that, for $j\in\{j_0-d_1,\ldots,j_0-1\}$ and $i$ such that $j \pmod{d_1} \equiv i$,
    \begin{align*}
        \Big| (\beta_i-\widetilde{\beta}_i) \ g^{(i)}_{\widetilde{\lambda}_j} \Big|
        &= \Big|\Big(c_{\widetilde{\lambda}_j} + \beta_i \ g^{(i)}_{\widetilde{\lambda}_j}\Big) - \Big(c_{\widetilde{\lambda}_j} + \widetilde{\beta}_i \ g^{(i)}_{\widetilde{\lambda}_j} \Big)\Big|\\
        &\leq \Big|c^{\beta}_{\widetilde{\lambda}_j}\Big| + \Big|c^{\widetilde{\beta}}_{\widetilde{\lambda}_j}\Big|
        \leq 2C 2^{-(j_0-d_1)\gamma}.
    \end{align*}
    So, by \eqref{eq:defi_sat_func:wlt_coef} and \eqref{eq:defi:sat_func_multi_dim_subspace}, $\big|\beta_i-\widetilde{\beta}_i\big| \leq 2C \frac{2^{-(j_0-d_1)\gamma}}{j^{-2/q}~\mu(\widetilde{\lambda}_j)}$.
    As for every $j\in\{j_0-d_1,\ldots,j_0-1\}$, one has $j>M_2$, then $\mu\big(\widetilde{\lambda}_j\big) \geq 2^{-j(\gamma-1/p)}$, and we deduce that $\mu\big(\widetilde{\lambda}_j\big) \geq 2^{-j_0(\gamma-1/p)}$.
    Going through all $j\in\{j_0-d_1,\ldots,j_0-1\}$ induce that for every $i\in \{1,\ldots,d_1\}$, one has
    \begin{align*}
        \big|\beta_i-\widetilde{\beta}_i\big|
        &\leq 2C \frac{2^{-(j_0-d_1)\gamma}}{j^{-2/q}~2^{-j_0(\gamma-1/p)}} 
        \leq 2C j_0^{2/q}~\frac{2^{-(j_0-d_1)\gamma}}{2^{-j_0(\gamma-1/p)}} \\
        &\leq 2C 2^{d_1\gamma}~2^{-j_0 \big(\frac{1}{p}-\frac{2}{q} \frac{\log_2(j_0)}{j_0}\big)}.
    \end{align*}
    Considering $j_0$ large enough, for every $i\in \{1,\ldots,d_1\}$, $\big|\beta_i-\widetilde{\beta}_i\big| \leq C 2^{-j_0/2p}$.
    This means that the distance between two elements of $\mathcal{B}_{\lambda_0}$ is at most $2^{-j_0/2p}$ up to a universal constant.
    Hence, the volume of $\mathcal{B}_{f,\lambda_0}$ is given by $\mathcal{L}_{d_1}(\mathcal{B}_{\lambda_0}) \leq C\ 2^{-j_0 \cdot d_1/2p}$.
    Summing over all the $2^{d j_0}$ dyadic cubes for $k_0 \in \{0,\ldots,2^{j_0}-1\}^d$ at scale $j_0$ gives $\mathcal{L}_{d_1} \big(\bigcup_{k_0 \in \Z^d_{j_0}} \mathcal{B}_{\lambda_0}\big) \leq C 2^{j_0(d-d_1/2p)}$.
    By our choice \eqref{eq:def_d_1}, $d<d_1/2p$ so the Borel-Cantelli lemma implies that 
    \begin{equation*}
        \mathcal{L}_{d_1} \bigg(\limsup_{j_0\to \infty} \bigcup_{k_0 \in \Z^d_{j_0}} \mathcal{B}_{\lambda_0}\bigg) = 0.
    \end{equation*}
    Finally by \eqref{eq:dec_B}, $\mathcal{L}_{d_1}(\mathcal{B}) = 0$, which is the claim of Proposition \ref{prop:beta_B_meas_0}.
\end{proof}

We can now conclude regarding the shyness of the sets $\mathcal{O}^{p,\gamma,M}$.
\begin{prop}
    For all $p\in\N^*$, $M\in(\N^{*})^2$ and $\gamma\in\Q_+$, $\mathcal{O}^{p,\gamma,M}$ is shy.
\end{prop}
\begin{proof}
    Let $\vartheta = G_{\#}\mathcal{L}_{d_1}$ be the push-forward measure of the $d_1$-dimensional Lebesgue measure by 
    \begin{equation}
        G:\bigg\{\begin{array}{ccc}
            [0,1]^{d_1} & \to & \mathcal{P} \\
            \beta & \mapsto & \sum_{i=1}^{d_1} \beta_i \cdot \mathscr{G}_{q}^{(i)}.
        \end{array}
    \end{equation}
    This measure is supported by a compact set of $\mathcal{P} \subset \widetilde{B}^{\mu,\Psi}_{\infty,q}([0,1]^d)$.
    
    Fix any $f\in\widetilde{B}^{\mu,\Psi}_{\infty,q}([0,1]^d)$. 
    Using Proposition \ref{prop:beta_B_meas_0}, for $\vartheta$-almost all $F\in\mathcal{P}$, $f+F \notin\mathcal{O}^{p,\gamma,M}$. 
    Hence for any $f\in\widetilde{B}^{\mu,\Psi}_{\infty,q}([0,1]^d)$, the set $\{f+\mathcal{O}^{p,\gamma,M}\}$ has $\vartheta$-measure equal to 0, i.e.,
    $\vartheta(\{f+\mathcal{O}^{p,\gamma,M}\})=0$.
    Since it is true for any $f\in \widetilde{B}^{\mu,\Psi}_{\infty,q}([0,1]^d)$, the set $\mathcal{O}^{p,\gamma,M}$ is shy.
\end{proof}

\noindent By Proposition \ref{prop:incl_cFxi_in_O}, $(\mathcal{F}^{\mu})^\complement$ is shy as subset of countable unions and intersec\-tions of shy sets. This yields that $\mathcal{F}^{\mu}$ is prevalent, hence Proposition \ref{prop:prev_set_besov_inf_q}. 

\section{Extensions}

We have identified the singularity spectrum for a prevalent set of functions in the case $p=\infty$. This leaves the finite case $p<\infty$, for which the computation of the singularity spectrum of functions requires, not only the singularity spectrum of the set function $\mu$, but also generalised results in Diophantine approximation and dynamical systems.
However, following Barral-Seuret's work \cite[Theorem 4]{Barral-Seuret:2023:Besov_Space_Part_2}, we conjecture an analogous result in the case of prevalence.
\begin{conj}
    Let $\mu$ be an almost doubling capacity on $[0,1]^d$ satisfying $(P)$ and verifying the SMF, let $p,q\in[1,+\infty]$.
    Consider the mappings
    \begin{equation*}
        \tau_{\mu,p}(t)=
        \begin{cases}
           \, \frac{p-t}{p}~\tau_{\mu}\Big( \frac{p}{p-t}~t \Big) &\text{if}~t\in(-\infty,p) \\
         \ \ d &\text{if}~t\in[p,+\infty)
        \end{cases}
    \end{equation*}
    and $\theta_p:\alpha\in[\tau^{\prime}_{\mu}(+\infty),\tau^{\prime}_{\mu}(-\infty)] \longmapsto \alpha + \frac{\tau^{*}_{\mu}(\alpha)}{p}$.
    For all $\displaystyle f\in \widetilde{B}^{\mu,\Psi}_{p,q}([0,1]^d)$,
    \begin{equation*}
        \sigma_f(h) \leq 
        \begin{cases}
           \, \sigma_{\mu}(\theta^{-1}(h)) = \tau_{\mu,p}^{*}(h) &\text{if}~h \leq \tau^{\prime}_{\mu,p}(0) \\
         \ \ d &\text{if}~h > \tau^{\prime}_{\mu,p}(0).
        \end{cases}
    \end{equation*}
  There is a prevalent set of $\widetilde{B}^{\mu,\Psi}_{p,q}([0,1]^d)$ satisfying $\sigma_f=\sigma_{\mu}(\theta^{-1}(\cdot))=\tau_{\mu,p}^{*}.$
\end{conj}

A second possible extension is to relax the hypothesis of Theorem \ref{theo:prev_spect_inh_besov_inf_q} to include capacities that are not almost doubling, such as the Bernoulli multinomial measures and multiplicative cascade.

\bibliographystyle{abbrv}
{
\setstretch{0}
\setlength{\bibsep}{5pt}
\bibliography{biblio.bib}
}

\end{document}